\documentclass[12pt]{article}

\usepackage[UKenglish]{babel}
\usepackage{amsmath,amsthm,amssymb,latexsym,mathrsfs}
\usepackage{longtable,epsfig,hhline,float,listings,color}
\usepackage{graphicx,caption,subcaption,bm,fancybox}
\usepackage[margin=1in]{geometry}

\usepackage[space]{grffile}
\graphicspath{{./Images/}}

\usepackage[breaklinks=true,colorlinks=true,linkcolor=blue,urlcolor=blue,citecolor=blue]{hyperref}

\newcommand{\red}[1]{\textcolor{black}{#1}}

\def\rd{\text{d}}
\def\re{\text{e}}
\def\ri{\text{i}}
\def\Zh{\widehat Z}
\def\eps{{\varepsilon}}

\def\qand{\quad\mbox{and}\quad}
\def\fr{\mbox{$\frac{1}{2}$}}
\def\lth{\langle\!\langle}
\def\rth{\rangle\!\rangle}
\def\blth{\left\langle\!\left\langle}
\def\brth{\right\rangle\!\right\rangle}

\def\R{\mathbb{R}}
\def\C{\mathbb{C}}
\def\T{\mathbb{T}}
\def\D{\mbox{${\rm D}$}}

\def\bk{{\bf k}}

\def\bw{{\bm \omega}}
\def\bq{{\bf q}}

\def\bO{{\bm \Omega}}
\def\be{{\bm \zeta}}
\def\bth{{\bm \theta}}
\def\bph{{\bm\phi}}
\def\bps{{\bm \psi}}
\def\bgam{{\bm \gamma}}

\def\Ups{\Upsilon}
\def\g{{\mathrm{g}}}

\begin{document}

\begin{center}
  \textsf{\Large Nonlinear theory for coalescing characteristics}\\[2mm]
  \textsf{\Large in multiphase Whitham modulation theory}
%   \textsf{\large generating nonlinearity, dispersion, and complexity}\\[2mm]
\vspace{0.75cm}

\textit{Thomas J. Bridges$\ ^1$ \& Daniel J. Ratliff$\ ^2$}
\vspace{0.25cm}

\textit{$\ ^1$Department of Mathematics, University of Surrey, Guildford GU2 7XH, UK}
\vspace{.1cm}

\textit{$\ ^2$Department of Mathematical Sciences, Loughborough University,
  Loughborough, Leicestershire LE11 3TU, UK}
\vspace{1.0cm}

\setlength{\fboxsep}{10pt}
\doublebox{\parbox{14cm}{
    {\bf Abstract.}
    The multiphase Whitham modulation equations with $N$ phases have
    $2N$ characteristics which may be of
    hyperbolic or elliptic type.  In this paper
    a nonlinear theory is developed for coalescence, where
    two characteristics change from hyperbolic to elliptic via collision.
    Firstly, a linear theory develops
    the structure of colliding characteristics involving the topological
    sign of characteristics and multiple Jordan chains,
and secondly a nonlinear modulation theory is developed for transitions.
The nonlinear
theory shows that coalescing characteristics morph the Whitham equations
into an asymptotically valid geometric form of the two-way Boussinesq equation.
That is, coalescing characteristics generate dispersion, nonlinearity
and complex wave fields. For illustration,
the theory is applied to coalescing characteristics associated with
the modulation of two-phase travelling-wave solutions
of coupled nonlinear Schr\"odinger equations,
highlighting how collisions can be identified and the relevant
 dispersive dynamics constructed.
}}
\end{center}
\vspace{.25cm}
\hfill{\footnotesize{\textsf{\today}}}
\vspace{-0.9cm}

%-------------------------------------
%{\footnotesize\tableofcontents}
%-------------------------------------

\section{Introduction}
\label{sec-intro}
\setcounter{equation}{0}

The theory of modulation, particularly Whitham modulation theory,
takes existing nonlinear waves, such as finite-amplitude periodic
travelling waves, and provides a framework for studying the dynamical
implications of perturbing the basic properties of the
nonlinear wave.  In classical
modulation, the properties of the basic state (wavenumber, frequency,
meanflow) are allowed to depend on space and time, and partial
differential equations (PDEs) are derived for these parameters.  Study
of these PDEs then provides information about the evolution of the
basic state under perturbation.

Given a basic state, there are several strategies for deriving 
modulation PDEs (averaging the Lagrangian, averaging conservation
laws, geometric optics ansatz, other ans\"atze).  In all cases the governing
equations produced by Whitham modulation theory (WMT), for a simple
one-phase periodic travelling wave, can be expressed in the canonical form
\begin{equation}\label{wmes}
q_T=\Omega_X \qand \frac{\partial\ }{\partial T}\mathscr{A}(\omega+\Omega,
k+q) + \frac{\partial\ }{\partial X}\mathscr{B}(\omega+\Omega,
k+q) =0\,.
\end{equation}
They are a pair of nonlinear first-order PDEs for the two unknowns
$\Omega(X,T)$, the modulation frequency, and $q(X,T)$, the modulation
wavenumber.
The parameters $(\omega,k)$ are representative of the wavetrain from which the
Whitham modulation equations are obtained, and $X=\eps x$ and
$T=\eps t$ are slow time and space scales.  The first equation
 is called \emph{conservation of waves} and the second
is called \emph{conservation of wave action} \cite{whitham-book}.
When the governing equations are the Euler-Lagrange equations associated
with a Lagrangian functional, the scalar-valued functions
$\mathscr{A}$ and $\mathscr{B}$ are related via
\begin{equation}\label{ABDefn}
\mathscr{A}=\mathscr{L}_\omega\,, \quad \mathscr{B}=\mathscr{L}_k\,.
\end{equation}
The function $\mathscr{L}(\omega,k)$ is obtained by averaging the Lagrangian
evaluated on the
periodic travelling wave with frequency $\omega$ and wavenumber $k$. 

The pair of quasilinear first-order equations (\ref{wmes}) can be 
classified based on their characteristics. The Whitham modulation
equations (WMEs) can either be
hyperbolic (real characteristics) or elliptic (complex characteristics)
and the transition signals a change of stability of the underlying
periodic waves \cite{whitham65,whitham-book,br17,br18}.  It is this change
of type, and its generalization to multiphase wavetrains, and its
nonlinear implications, that are the main themes of this paper.

To identify the structure of coalescing characteristics, first consider
the one-phase case where only two characteristics exist and so
coalescence is elementary.
The linearization of the one-phase WMEs (\ref{wmes})
about the basic state, represented by $(\omega,k)$, is
\begin{equation}\label{wmes-linear-intro}
  q_T = \Omega_X \qand \mathscr{A}_\omega \Omega_T +
  \mathscr{A}_k q_T + \mathscr{B}_\omega\Omega_X + \mathscr{B}_k q_X=0\,,
\end{equation}
or, under the assumption $\mathscr{A}_\omega\neq0$, they can be written
in the standard hydrodynamical form,
\begin{equation}\label{wme-3}
\begin{pmatrix} q\\ \Omega \end{pmatrix}_T + {\bf F}(\omega,k)
\begin{pmatrix} q\\ \Omega \end{pmatrix}_X = \begin{pmatrix}0\\ 0\end{pmatrix}\,,
\end{equation}
with
\begin{equation}\label{M-def-1}
{\bf F}(\omega,k) = \frac{1}{\mathscr{A}_\omega}\left[\begin{matrix}
0 & - \mathscr{A}_\omega \\ \mathscr{B}_k & \mathscr{A}_k+\mathscr{B}_\omega\end{matrix}
\right]\,.
\end{equation}
Here, $\mathscr{A}$ and $\mathscr{B}$ are evaluated at $\Omega=q=0$.
The characteristics (eigenvalues of ${\bf F}$) are
\begin{equation}\label{wmes-characteristics}
c^{\pm} = \frac{\mathscr{A}_k+\mathscr{B}_\omega}{2\mathscr{A}_\omega}
\pm \frac{1}{\mathscr{A}_\omega}\sqrt{-\Delta_L}\,,
\end{equation}
where
\begin{equation}\label{Delta-L-def}
  \Delta_L = \mathscr{A}_\omega\mathscr{B}_k-\mathscr{A}_k\mathscr{B}_\omega=
{\rm det}\left[\begin{matrix} \mathscr{L}_{\omega\omega} & \mathscr{L}_{\omega k}\\
\mathscr{L}_{k\omega} & \mathscr{L}_{kk}\end{matrix}\right]\,,
\end{equation}
using (\ref{ABDefn}) in the latter equality.  The sign of the
determinant $\Delta_L$,
called the \emph{Lighthill determinant} (\textsc{Lighthill}~\cite{lighthill67}),
signals whether the characteristics are real or complex,
\[
\begin{array}{rcl}
  \Delta_L<0 \quad &\Longrightarrow& \quad \mbox{hyperbolic WMEs}\\[2mm]
   \Delta_L>0 \quad &\Longrightarrow& \quad \mbox{elliptic WMEs}\,.
\end{array}
\]
At the transition, when $\Delta_L=0$, the two characteristics are
equal, \red{the} Whitham modulation \red{equations degenerate}, and a new modulation
strategy is needed.
In \cite{br17} a nonlinear modulation theory is developed for
\red{the above case within} the WMEs in the case of one-phase wavetrains.  It is
valid near the transition from hyperbolic to elliptic, showing
that the WMEs (\ref{wmes}) are replaced by
\begin{equation}\label{2way-boussinesq-lps}
q_T=\Omega_X \qand
\mathscr{A}_\omega \Omega_T + \kappa qq_X + \mathscr{K}q_{XXX}=0\,,
\end{equation}
where $T=\eps^2t$, $X=\eps(x-c_gt)$, and $c_g$ is a nonlinear
group velocity at the transition.  The coefficients $\mathscr{A}_\omega$
and $\kappa$ are obtained from derivatives of the components of
conservation of wave action, and the dispersion coefficient
$\mathscr{K}$ arises due to a symplectic Jordan chain argument.
Differentiating the second equation of (\ref{2way-boussinesq-lps})
with respect to $X$ and using
the first equation reveals that it
is a variant of the two-way Boussinesq equation for $q$,
\begin{equation}\label{q-boussinesq-intro}
\mathscr{A}_\omega\, q_{TT}+ \left(\fr\kappa q^2 + \mathscr{K}\, q_{XX}\right)_{XX} = 0\,.
\end{equation}
The coefficients in (\ref{2way-boussinesq-lps}) and
(\ref{q-boussinesq-intro}) are universal in the same sense that the Whitham
equations are universal  -- they follow from abstract properties of
a Lagrangian.  Extension of the derivation of
(\ref{2way-boussinesq-lps}) to two space dimensions and time appears
in \cite{br18}. The emergence of the equation
(\ref{q-boussinesq-intro}) shows that
coalescing characteristics generate nonlinearity, dispersion and wave
fields of greater complexity. The complexity is due to the
wide range of known localized, multi-pulse, quasiperiodic, and extreme value
solutions of the two-way Boussinesq equation.

In order \red{to} generalize this nonlinear theory for coalescing characteristics
to the case of multiphase wavetrains several new results are needed.
\red{
  The first results on non-generic Whitham modulation theory, in
  the multiphase case, considered the case when
  the generic WMEs have a single or double zero characteristic.
In \textsc{Ratliff \& Bridges}~\cite{rb16a},
it was shown that a zero characteristic in WMT leads, under
re-modulation, to Korteweg-de Vries (KdV) dynamics on a longer
time scale.  In \textsc{Ratliff}~\cite{r18}, 
a double zero characteristic is re-modulated leading to dynamics governed
by a two-way Boussinesq equation.}

\red{The motivation for re-modulation, in both cases,
is that a zero characteristic in the generic WMEs suggests no dynamics,
but in fact the dynamics is moved to a slower time scale.  The time
scale $t\thicksim \eps^{-1}$, in generic WMT,
is replaced by $t\thicksim\eps^{-2}$ in \cite{r18}
and is replaced by $t\thicksim\eps^{-3}$ in \cite{rb16a}.
The double zero characteristic is a special case of coalescing characteristics.
It does not require the theory of 
coalescing characteristics with non-zero speeds
such as the sign characteristic, and it is
codimension two.}

\red{In this paper we are interested in the nonlinear theory
  near coalescing characteristics
  with non-zero speed.
This case is codimension one and so more likely to occur in applications,
and it requires the theory of the sign characteristic to track
collisions of characteristics. It was discovered in \cite{br19} that
every characteristic in the Whitham theory carries a topological sign,
and this sign is an important diagnostic as only coalescing
characteristics with opposite sign can change type from hyperbolic
to elliptic.}
In addition, several facets of the linear theory, such as intertwining Jordan
chains, that generate the coefficient $\mathscr{K}$, bring in new challenges.
For the nonlinear theory, we find that the form of the two-way Boussinesq
equation (\ref{q-boussinesq-intro}) carries over to the case of
coalescing characteristics with nonzero speed, 
but there is a discrepancy between the fact that (\ref{q-boussinesq-intro}) is
scalar valued but the WMEs in the multiphase case have $2N$ equations.
Hence a secondary reduction of the nonlinear equations is required.
Showing that the coefficients are universal is also an order of magnitude
more difficult in this case.

\begin{figure}[ht]
\begin{center}
\includegraphics[width=6.0cm]{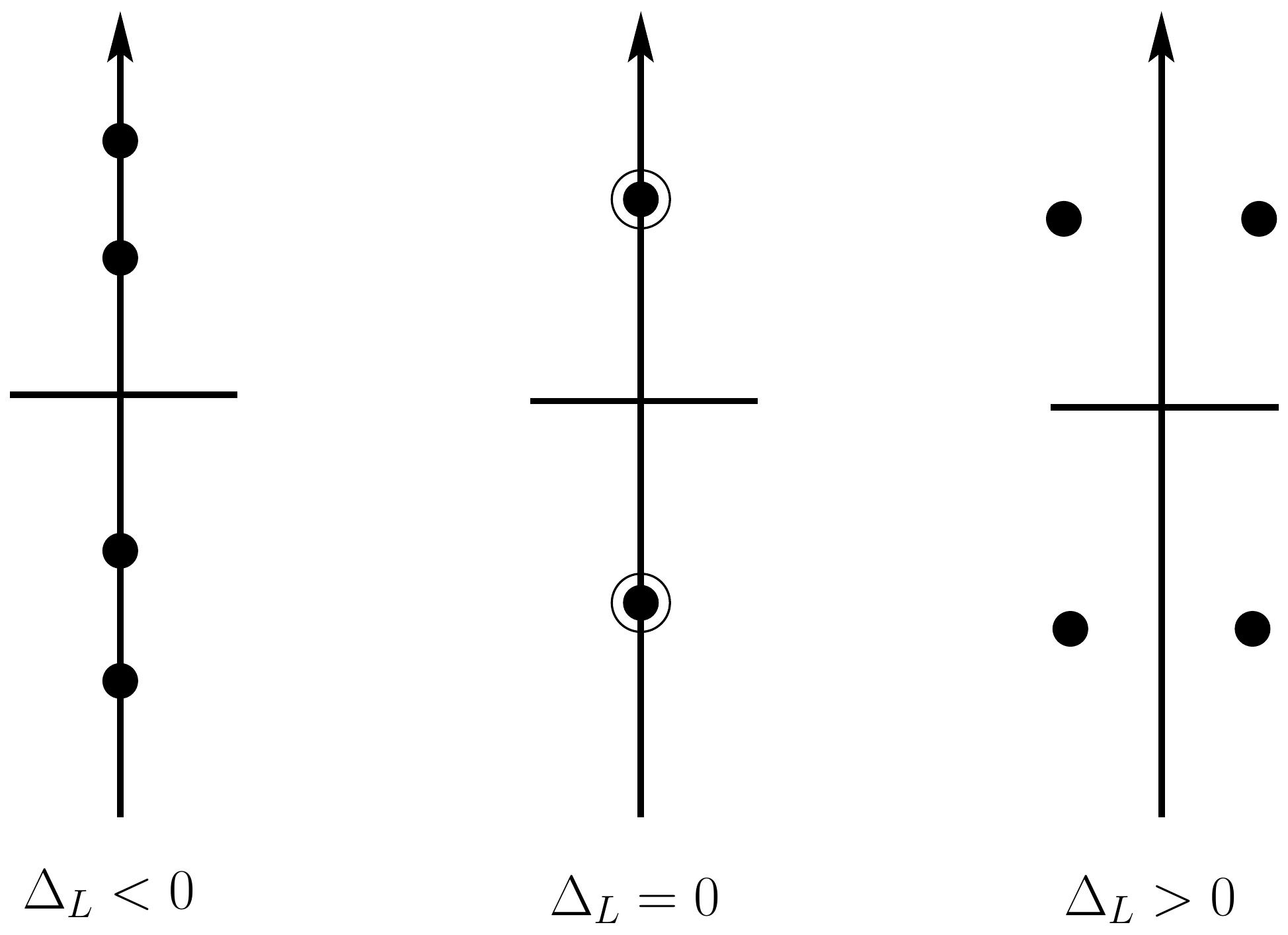}
\end{center}
\caption{Collision of purely imaginary eigenvalues in the Whitham equations.}
\label{fig_Delta-L-eigenvalues}
\end{figure}

The mathematics of how characteristics coalesce and change type
is addressed as follows.  Firstly consider
the one-phase case.
The change of type of the characteristics signals an instability of
the basic state, and this linear instability is made apparent by
taking the normal-mode ansatz
\begin{equation}\label{normal-modes}
\begin{pmatrix} q(X,T)\\ \Omega(X,T)\end{pmatrix} =
{\rm Re}\left\{ \begin{pmatrix} \widehat q\\ \widehat\Omega\end{pmatrix}
    \re^{\lambda T + \ri\nu X}\right\}\,,
\end{equation}
and substituting into (\ref{wmes-linear-intro}) to obtain
\begin{equation}\label{lambda-eigs}
       \lambda = \red{\pm}\ri c^{\pm}\nu \red{=
\pm \ri \nu\left( \frac{\mathscr{A}_k+\mathscr{B}_\omega}{2\mathscr{A}_\omega}
\pm \frac{1}{\mathscr{A}_\omega}\sqrt{-\Delta_L}\right)}\,.
\end{equation}
\red{There are four $\lambda-$eigenvalues for fixed $\nu\neq0$ since
$\widehat q$ and $\widehat\Omega$ are complex valued.}
An unstable exponent (${\rm Re}(\lambda)>0$)
with modulation wave number $\nu$ exists precisely when $\Delta_L>0$.
As $\Delta_L$ changes sign, the eigenvalues (\ref{lambda-eigs})
change from four purely-imaginary eigenvalues
to a complex quartet as shown schematically in Figure
\ref{fig_Delta-L-eigenvalues}.
This type of stability transition is familiar from the theory of linear
Hamiltonian systems, as it is precisely the Hamiltonian Hopf bifurcation~\cite{v85},
 and in that setting the collision and resulting instability
occurs since the eigenvalues have opposite Krein signature \cite{howard,v85}.
However, as shown in \cite{br19}, there is no obvious symplectic structure
in the Whitham theory,
and it is the \emph{sign characteristic of Hermitian matrix pencils}
that is operational here.  The sign characteristic has a central
role in the theory of Hermitian matrix pencils relative to an indefinite
metric (see \textsc{Gohberg et al.}~\cite{glr-book} for a history and
references).

The Hermitian matrix pencil structure of (\ref{wmes-linear-intro})
is evoked by multiplying the conservation of waves by $\mathscr{A}_\omega$,
assuming $\mathscr{A}_\omega \neq 0$, and combining the two
equations in (\ref{wmes-linear-intro}) as
\begin{equation}\label{wme-3-a}
  \left[\begin{matrix} 0 & \mathscr{A}_\omega \\
      \mathscr{A}_\omega & \mathscr{A}_k+\mathscr{B}_\omega
    \end{matrix}\right] \begin{pmatrix} \Omega\\ q \end{pmatrix}_T
  + \left[\begin{matrix} -\mathscr{A}_\omega & 0 \\ 0 & \mathscr{B}_k
      \end{matrix}\right]
  \begin{pmatrix} \Omega\\ q \end{pmatrix}_X
  = \begin{pmatrix}0\\ 0\end{pmatrix}\,.
\end{equation}
The two coefficient matrices are symmetric.  Now the modified normal mode
ansatz
\[
\begin{pmatrix} \Omega(X,T) \\ q(X,T) \end{pmatrix} =
{\rm Re}\left\{ \begin{pmatrix} \widehat\Omega\\ \widehat q \end{pmatrix}
\re^{\ri\nu(X+cT)}\right\}\,,
\]
generates the following Hermitian matrix eigenvalue problem
\begin{equation}\label{wme-3-b}
\left( \left[\begin{matrix} -\mathscr{A}_\omega & 0 \\ 0 & \mathscr{B}_k
      \end{matrix}\right] +c
\left[\begin{matrix} 0 & \mathscr{A}_\omega \\
      \mathscr{A}_\omega & \mathscr{A}_k+\mathscr{B}_\omega
    \end{matrix}\right] \right)\begin{pmatrix} \widehat\Omega\\ \widehat q \end{pmatrix}  = \begin{pmatrix}0\\ 0\end{pmatrix}\,.
\end{equation}
The theory of Hermitian matrix pencils shows
that each eigenvalue of (\ref{wme-3-b})
has a \emph{sign characteristic} and a necessary condition
for instability is that eigenvalues coalesce and have opposite sign
characteristic \cite{br19}.  In the one phase case,
with just two characteristics, the sign characteristic is less
interesting, and indeed trivial.  In the multiphase case,
with many characteristics, the
coalescence of characteristics may or may not lead to instability,
and so the sign characteristic becomes an essential diagnostic tool.
The principal case of
interest in this paper is when all the characteristics are real,
with only one pair, having opposite sign,
undergoing a transition to instability.

The generalization of the dispersionless WMEs (\ref{wmes}) to the
multiphase case is
\begin{equation}\label{mphase-wmt-canonical}
  \bq_T = \bO_X \qand
  \frac{\partial\ }{\partial T}\mathbf{A}(\bw+\bO,
\bk+\bq) + \frac{\partial\ }{\partial X}\mathbf{B}(\bw+\bO,\bk+\bq) =0\,,
\end{equation}
where $\bw,\,\bk\in\R^N$ are given parameters representative
of the basic state, and $\bq,\bO\in\R^N$ are the vector-valued unknowns,
the modulation wavenumber and frequency, which depend on $T=\eps t$ and $X=\eps x$. 
\red{This general form of the WMEs (\ref{mphase-wmt-canonical})
  covers not only multiply-periodic waves, but quasi-periodic structures,
  and WMEs on pseudo-phases which model meanflow
  (this latter case is discussed in \S\ref{subsec-pseudo}).}
When the governing equations are the Euler-Lagrange equation associated
with a Lagrangian functional, the mappings ${\bf A}$ and ${\bf B}$ are \red{again}
variation\red{s of the average Lagrangian $\mathscr{L}(\bw,\bk)$} with the properties,
\begin{equation}\label{A-Lw-def}
\mathbf{A}(\bw+\bO,\bk+\bq) =
\D_{\bw}\mathscr{L}(\bw+\bO,\bk+\bq)\,,
\end{equation}
and
\begin{equation}\label{B-Lk-def}
  \mathbf{B}(\bw+\bO,\bk+\bq) =
\D_{\bk}\mathscr{L}(\bw+\bO,\bk+\bq)\,.
\end{equation}
Cross differentiating shows that the Jacobians satisfy
\begin{equation}\label{AB-symmetry}
  \D_\bk{\bf A} = \big( \D_\bw{\bf B}\big)^T\,.
\end{equation}
This symmetry will be important for generalising the Hermitian
property of (\ref{wme-3-b}) to the multiphase case.

Given a smooth \red{averaged Lagrangian} $\mathscr{L}$,
the pair of equations (\ref{mphase-wmt-canonical}) is a closed
first-order system of PDEs for $\bO$ and $\bq$ with \red{up to} $2N$ characteristics.
This formulation of vector-valued WMEs was introduced
in \textsc{Ratliff}~\cite{r18,r17b}.  However, multiphase Whitham modulation
theory has a rich history. Multiphase WMEs
were first introduced and studied by
\textsc{Ablowitz \& Benney}~\cite{ab70}, in the context of
scalar nonlinear wave equations, where the appearance of small
divisors was noted.
For integrable systems small divisors disappear:
multiphase averaging and the Whitham equations are robust and rigorous,
and a general theory can be obtained
(e.g.\ \textsc{Flashka et al.}~\cite{ffm80} and its citation trail).
\textsc{Whitham}~\cite[\S 14.7]{whitham-book}
includes potential variables as additional phases (``pseudo-phases''
\red{which
are included in the theory here; see \S\ref{subsec-pseudo}})
and generates a form of multiphase modulation and applies it to
wave-meanflow interaction of Stokes water waves \cite{whitham67}.
\textsc{Willebrand}~\cite{willebrand} takes
multiphase modulation theory to a new level by deriving the $N-$phase
WMEs, for Stokes wave solutions of the water wave problem, with $N$ arbitrary
and he takes the limit $N\to\infty$. This theory is formal
and the series are divergent and have small divisors, but the
leading order terms are instructive (see comments on this later in
\S\ref{sec-infinite-phase}).
The theory of \cite{willebrand} is now used in ocean wave forecasting
  (e.g.\ Chapter
  9 of \textsc{Olbers et al.}~\cite{owe12}).

On the other hand, when the system is not integrable, but there is an $N-$fold
symmetry, a theory for conservation of wave action can be
developed without small divisors and smoothly varying
$N-$phase wavetrains. This strategy is implemented in
\cite{r18,r17b}, where multiphase
wavetrains are characterized as relative equilibria with
smooth dependence on parameters.

Going back to the abstract multiphase WMEs (\ref{mphase-wmt-canonical}),
with the symmetry property (\ref{AB-symmetry}) and the gradient properties
(\ref{A-Lw-def}), the linearization of (\ref{mphase-wmt-canonical})
can be cast into the form of a Hermitian matrix pencil,
\begin{equation}\label{wme-3-c}
\left[ \left[\begin{matrix} -\D_\bw{\bf A} & 0 \\ 0 & \D_{\bk}{\bf B}
      \end{matrix}\right] +c
\left[\begin{matrix} 0 & \D_\bw{\bf A} \\
      \D_\bw{\bf A} & \D_\bk{\bf A}+\D_\bw{\bf B}
    \end{matrix}\right] \right]\begin{pmatrix} \widehat\bO\\ \widehat \bq \end{pmatrix}  = \begin{pmatrix}{\bf 0}\\ {\bf 0}\end{pmatrix}\,,
\end{equation}
assuming that $\D_\bw{\bf A}$ is invertible.  The Jacobians
$\D_\bw{\bf A}$, \red{$\D_\bk{\bf B}$}, $\D_\bw{\bf B}$, and \red{$\D_\bk{\bf A}$},
are $N\times N$ matrices, with the first two symmetric and the latter two
\red{are the matrix transpose of one another}.
The $2N\times 2N$ linear eigenvalue problem (\ref{wme-3-c}) can be reduced,
by eliminating $\widehat\bO$,
\begin{equation}\label{Omega-q-def}
  \widehat\bO = c\widehat\bq\,,\quad(\mbox{assuming}\ {\rm det}[\D_\bw{\bf A}]\neq0)\,,    
\end{equation}
to an $N\times N$ \emph{quadratic Hermitian matrix pencil},
\begin{equation}\label{wme-3-quadratic}
  {\bf E}(c)\widehat\bq :=
  \left[ \D_\bw{\bf A}\, c^2 + c(\D_\bk{\bf A}+\D_\bw{\bf B})
    + \D_\bk{\bf B}\right] \widehat \bq  = {\bf 0}\,.
\end{equation}
A parallel theory can be developed for the
sign characteristic in this context \cite{glr80,mntx16,tm01}.
Suppose $c_0$ is a simple real eigenvalue satisfying ${\rm det}[{\bf E}(c_0)]=0$ with
eigenvector \red{$\be$, so that}
\begin{equation}\label{zeta-0-def}
  {\bf E}(c_0)\be = 0 \,.
\end{equation}
Then the sign characteristic of $c_0$ is
\begin{equation}\label{sign-char-def}
S(c_0) = {\rm sign}\left( \left\langle \be,{\bf E}'(c_0)\be\right\rangle\right)\,,
\end{equation}
where $\langle\cdot,\cdot\rangle$ is an inner product on $\R^N$,
and the prime denotes differentiation with respect to $c$.
A discussion of the history and various formulations of the sign characteristic
is given in \cite{br19}. 
The sign characteristic is invariant under congruence transformation,
        ${\bf E}(c_0)\mapsto {\bf P}^T{\bf E}(c_0){\bf P}$, for
        any invertible ${\bf P}$ \cite{glr80}.
The quadratic formulation (\ref{wme-3-quadratic}), rather
than its linearisation (\ref{wme-3-c}), turns out to be the most efficient
in applications and arises naturally in the modulation theory.
        
Starting with the quadratic Hermitian matrix pencil (\ref{wme-3-quadratic})
a theory for the sign characteristic in the context of Whitham modulation theory
is developed by \textsc{Bridges \& Ratliff}~\cite{br19}.
The $2N$ characteristics of the linearized problem satisfy
\begin{equation}\label{Ec-def}
\Delta(c) := {\rm det}[ {\bf E}(c)]=0\,.
\end{equation}
Double non-semisimple characteristics, which characterize
coalescence, and define the nonlinear group velocity $c_g$, satisfy
\begin{equation}\label{cg-def-intro}
\Delta(c_g)=\Delta'(c_g)=0 \quad\mbox{but}\quad \Delta''(c_g)\neq 0 \,,
\end{equation}
with a single geometric eigenvector
\begin{equation}\label{zeta-def}
  {\bf E}(c_g)\be = 0 \,,
\end{equation}
and generalized eigenvector
\begin{equation}\label{gamma-def}
  {\bf E}(c_g)\bgam = -{\bf E}'(c_g)\be\,.
\end{equation}
Solvability of (\ref{gamma-def}) defines $c_g$.  All these properties
follow from the structure of the linear operator ${\bf E}(c)$
with $c\in\R$ and are studied in \cite{br19} and the details required
here are recorded in \S\ref{sec-characteristics} and their role
in the Jordan chain theory is developed
in \S\ref{sec-jordanchain}.

When multiphase modulation is introduced for the nonlinear problem,
with an appropriate scaling, the vector-valued conservation of wave action
(\ref{mphase-wmt-canonical}) will be morphed into another form with
dispersion.  It will however still have dimension $N$, so a further reduction
is necessary in order to obtain a generalisation of (\ref{q-boussinesq-intro}).
The strategy is to split $\R^N={\rm span}\{\be\}\oplus \R^{N-1}$.
The geometric eigenvector $\be$, defined in (\ref{zeta-0-def}),
provides a preferred direction in $\bq-$wavenumber
space associated with the coalescence.
This preferred direction is an essential part of the nonlinear modulation theory.
It provides a projection operator so that the vector-valued
conservation of wave action (\ref{mphase-wmt-canonical}) can
be reduced to a scalar equation, and this scalar equation,
which also requires a rescaling of the slow variables, and extension of the analysis to fifth order in $\eps$, is a geometric form of the 
scalar-valued two-way Boussinesq equation
\begin{equation}\label{2-way-boussinesq-intro}
\mu\, U_{TT}+ \fr\kappa(U^2)_{XX} +\mathscr{K}\, U_{XXXX} = 0\,,
\end{equation}
where $\mu$ and $\kappa$ are determined by the geometry of
the averaged Lagrangian $\mathscr{L}(\bw,\bk)$ and $\mathscr{K}$ is
determined by a twisted Jordan chain argument. 

The geometry of
$\mathscr{L}(\bw,\bk)$ is discussed in \S\ref{sec-lagrangian-geometry}.  The most
remarkable outcome of the geometry is that the coefficient $\kappa$
in (\ref{2-way-boussinesq-intro}) has the simple formula
\begin{equation}\label{kappa-def}
  \kappa := \frac{d^3\ }{ds^3}\mathscr{L}(\bw + s c_g \be,\bk+s\be)\bigg|_{s=0}\,.
\end{equation}
  The coefficient $\mu$ is determined by a Jordan chain associated with the
  linear operator
  ${\bf E}(c_g)$ in (\ref{zeta-def}).  Indeed, $\mu\neq0$ is the
  condition required for termination of the Jordan chain $(\be,\bgam)$ in (\ref{zeta-def})-(\ref{gamma-def}).
  A different Jordan chain, associated
  with the linearization of the Euler-Lagrange equations (denoted ${\bf L}$
  in (\ref{L-def})), determines the dispersion coefficient $\mathscr{K}$.
  This latter Jordan chain argument is similar to the case of
  multiphase modulation associated with zero characteristics in \cite{r18,r17b,rb16a} but
here the Jordan chain intertwines two different chains generated by ${\bf L}$.
The nonlinear modulation theory where (\ref{kappa-def}) arises
naturally and feeds into the emergence of (\ref{2-way-boussinesq-intro})
is developed in \S\ref{sec-mod-ansatz}.

The starting point for the theory is a general class of nonlinear PDEs generated
by a Lagrangian, and this class is introduced in \S\ref{sec-lagr-goveqns}.
Given a basic multiphase wavetrain $\widehat Z(\bth,\bw,\bk)$,
with phase $\bth=\bk x + \bw t + \bth_0$ and vector-valued frequency and
wavenumber $\bw,\bk$, satisfying the
Euler-Lagrange equations,
the dispersionless vector-valued WMEs (\ref{mphase-wmt-canonical})
are derived by modulating the basic state with
a geometric optics scaling \cite{r18,r17b,br19}.
The appropriate modulation ansatz is
\begin{equation}\label{geom-optics-ansatz}
Z(x,t) = \Zh(\bth +\eps^{-1}\bph,\bw + \bO,\bk+\bq) + \eps W(\bth +\eps^{-1}\bph,X,T,\eps)\,,
\end{equation}
where $\bph$, $\bw$ and $\bq$, depending on $T=\eps t$ and $X=\eps x$,
are the modulated phase, frequency and wavenumber, and $W$ is a remainder. 
Substitution of (\ref{geom-optics-ansatz})
into the Euler-Lagrange equation and
solvability requires $\bq$ and $\bO$ to satisfy (\ref{mphase-wmt-canonical})
to leading order \cite{r18,r17b}.

When two characteristics, of opposite sign characteristic, coalesce
and transition to instability, the geometric optics modulation ansatz
\[
\bth\mapsto\bth +\eps^{-1}\bph\,,\quad \bw\mapsto \bw + \bO\,,\quad \bk\mapsto
\bk+\bq\,,
\]
(with $T=\eps t$ and $X=\eps x$) in (\ref{geom-optics-ansatz}) must
be replaced.  The altered form utilised is
\begin{equation}\label{cc-ansatz}
  \bth \mapsto \bth+{\eps}{\bm\Phi}\,,\quad
  \bk \mapsto \bk+{\eps^2}{\bm\Phi}_X\,,\quad
  \bw \mapsto \bw+\eps^2c_g{\bm\Phi}_X + \eps^3{\bm\Phi}_T\,,
\end{equation}
where ${\bm\Phi}$ is a function of
the slow time and space variables,
\begin{equation}\label{slow-scales-intro}
X = {\eps}(x+c_gt)\,, \quad T = {\eps^2}t\,,
\end{equation}
with $c_g$ determined as part of the analysis,
and $\eps$ measuring the distance in $(\bw,\bk)-$space from the singularity (\ref{cg-def-intro}).  The new ansatz at coalescence is
\begin{equation}\label{Z-cc-ansatz}
Z(x,t) = \Zh(\bth +{\eps}{\bm\Phi},\bw +\eps^2c_g{\bm\Phi}_X+\eps^3{\bm\Phi}_T,\bk+\eps^2{\bm\Phi}_X) + \eps^3 W(\bth,X,T,\eps)\,.
\end{equation}
Finer detail on the ansatz, including definitions of $\bq$ and $\bO$ and their
relation to ${\bm\Phi}_X$ and ${\bm\Phi}_T$ is given in \S\ref{sec-mod-ansatz}.
Substitution of this ansatz into the governing Euler-Lagrange equations,
expanding everything in powers of $\eps$, and setting order by order to zero, results,
by imposing a solvability condition,
in a vector-valued two-way Boussinesq equation induced
by conservation of wave action.
The projection
operator, defined using ${\rm Ker}({\bf E}(c_g))$, is then implemented
to split the conservation of wave action into two parts, one
generating the two-way Boussinesq equation
(\ref{2-way-boussinesq-intro}) with the complementary
part carrying over to higher order.

The paper has four parts:
the Lagrangian (geometry, analysis, and Euler-Lagrange equation),
the linear theory (for both operators ${\bf E}(c)$ and ${\bf L}$),
the nonlinear modulation analysis (implementing the ansatz
(\ref{Z-cc-ansatz})), and an illustrative example.

In \S\ref{sec-lagr-goveqns} a class of Lagrangian functionals
and a class of basic states is introduced.
In \S\ref{sec-lagrangian-geometry} the geometry of the mapping
$(\bw,\bk)\mapsto\mathscr{L}(\bw,\bk)$, where
$\mathscr{L}(\bw,\bk)$ is the Lagrangian evaluated on a
basic state, is studied.  Remarkably, many
of the features of the linear problem as well as the nonlinear modulation
are determined by the geometry of this scalar-valued function.
We end this discussion by reviewing the Whitham modulation theory
from a geometric perspective in \S\ref{sec-generic-multiphase}
to demonstrate how the characteristics
and their coalescence may be formulated using these notions, as discussed in 
\S\ref{sec-characteristics}.  In \S\ref{subsec-pseudo} we show how
the theory of pseudo-phases, and their relation with meanflow,
fits into the theory of this paper.

The linear theory has two parts: the structure of
the linear operator ${\bf E}(c)$ in (\ref{wme-3-quadratic}), including
the Jordan chain theory in the setting of the quadratic eigenvalue
problem (\ref{wme-3-quadratic}).
This theory is developed in \S\ref{sec-lagrangian-geometry} and
\S\ref{sec-jordanchain},
and appeals to the theory of sign characteristic for Hermitian
matrix pencils developed in \cite{br19}.  The second part of the
linear theory is the linearization of the Euler-Lagrange equation, which
is needed to develop a secondary Jordan chain needed for constructing
the dispersion coefficient $\mathscr{K}$ and the nonlinear modulation
theory, and this theory is developed in \S\ref{sec-jordanchain}.

The nonlinear theory
is developed in \S\ref{sec-mod-ansatz}.  Although the ansatz
(\ref{Z-cc-ansatz}) with (\ref{cc-ansatz}) is new, once the ansatz is
identified the strategy is similar to our previous papers, particularly
\cite{br17} and \cite{rb16a}, and so only the key new features are highlighted.
The theory is illustrated by application to the two-phase travelling
wave solutions of a class of coupled nonlinear Schr\"odinger
(CNLS) equations.  In \cite{br19} it was shown that these travelling wave
solutions have coalescing characteristics with transition to instability.  Here
the theory is applied to show the emergence of a geometric two-way Boussinesq
equation at these singularities.  Potential generalizations are discussed
in the concluding remarks section.

\section{The Lagrangian and governing equations}
\label{sec-lagr-goveqns}
\setcounter{equation}{0}

The theory is built on a general class of Lagrangian functionals
\[
\mathcal{L}(V) = \int_{t_1}^{t_2}\int_{x_1}^{x_2} L(V,V_t,V_x,\ldots)\,\rd x\rd t\,,
\]
where $V(x,t)$ is a vector-valued smooth field defined on the
rectangle $[x_1,x_2]\times[t_1,t_2]$. The lower dots indicate that 
the Lagrangian may also depend on higher derivatives of $V$, and the subsequent
theory can be adapted for these cases.  \red{Potential variables
  are included by letting some components of $V$ appear in $L(\cdot)$
  with derivatives only (see discussion in \S\ref{subsec-pseudo}
  of pseudo-phases associated with potential variables).}

Normally a non-degeneracy condition
on derivatives of $L$ with respect to $V_t$ and $V_x$ is assumed, but here
these conditions are circumvented by assuming up front that the
Lagrangian has been transformed to standard multisymplectic form
\begin{equation}\label{Lagr-mss}
\mathcal{L}(Z) =
\int_{t_1}^{t_2}\int_{x_1}^{x_2} \left[ \fr \langle
    {\bf M}Z_t,Z\rangle +
    \fr \langle{\bf J}Z_x,Z\rangle- S(Z)\right] \,\rd x\rd t\,,
\end{equation}
        where $Z\in\R^n$, ${\bf M}$ and ${\bf J}$ are
        skew-symmetric matrices, $S:\R^n\to\R$ is a given smooth
        function, and $\langle\cdot,\cdot\rangle$ is a standard inner
        product on $\R^n$.  For definiteness, $n$ is taken to be even
        and
        \begin{equation}\label{JcM}
          {\rm det} \big[  {\bf J} + c{\bf M} \big]\neq 0\,,\quad \forall\ c\in\mathscr{C}\subset\R\,,
        \end{equation}
\red{where $\mathscr{C}$ is some open set of real numbers.}
          Examples with $n=4$ for ${\bf M}$ and ${\bf J}$
          include the dispersive shallow water equations~\cite{tjb-book,b13},
          \[
          {\bf M} = \begin{pmatrix}
          0&-1&0&0\\
          1&0&0&0\\
          0&0&0&0\\
          0&0&0&0
          \end{pmatrix}\,, \qand
          {\bf J} = 
          \begin{pmatrix}
          0&0&-1&0\\
          0&0&0&-1\\
          1&0&0&0\\
          0&1&0&0
          \end{pmatrix}\,,
          \]
\red{which has $\mathscr{C}=\R$},
          and the coupled-mode equation (also massive-Thirring equation),
          \[
  {\bf M} = \begin{pmatrix} 0 & -1 & 0 & 0 \\ 1 & 0 & 0 & 0 \\
      0 & 0 & 0 & 1 \\ 0 & 0 & -1 & 0 \end{pmatrix}\qand
   {\bf J} = \begin{pmatrix} 0 & 0 & 1 & 0 \\ 0 & 0 & 0 & -1 \\
       -1 & 0 & 0 & 0 \\ 0 & 1 & 0 & 0 \end{pmatrix}\,.
   \]
   with an appropriate  choice of $S(Z)$ in both cases.
   \red{The latter case has $\mathscr{C}=\R\setminus\{\pm 1\}$,
since ${\rm det}[{\bf J}+c{\bf M}]=(c^2-1)^2$.}
   These examples and others can be found in
\cite{bd01,bhl10,b13,tjb-book}.

The Euler-Lagrange equation associated with (\ref{Lagr-mss}) is
\begin{equation}\label{MJS}
{\bf M}Z_t+{\bf J}Z_x = \nabla S(Z)\,,\quad Z\in\R^n\,.
\end{equation}
The theoretical developments to follow are based on this abstract form of
the Euler-Lagrange equation, with ${\bf M},{\bf J}$
general skew-symmetric matrices satisfying (\ref{JcM}),
\red{and $n$ even with $n\geq 2N$.}

\subsection{Symmetry, relative equilibria, and the basic state}
\label{subsec-symmetry}

The easiest way to generate smooth families of multiphase wavetrains is
to consider a Lagrangian that is invariant under the action of a Lie group.
Here, and henceforth it is assumed that the Lie group is abelian,
and some combination of $S^1$ (associated with periodic wave trains),
and affine translations.
\red{Affine translations, which are associated with pseudo-phases,
  are captured in
  the formulation (\ref{MJS}) by
  \[
  \langle {\bf e}_j,\nabla S(Z)\rangle =0\,,\quad j=1,\ldots,P\,,
  \]
  where $P$ is the dimension of the
  affine symmetry group, and ${\bf e}_j\in\R^n$ are the generators.}
The theory will be developed for the case of the $N-$torus,
which is appropriate for periodic $N-$phase wavetrains, as
the affine translation group is
much simpler and the necessary changes will
be recorded when needed.

Assume that
(\ref{MJS}) is equivariant with respect to an $\red{N}-$torus,
$\T^{\red{N}}=S^1\times S^1$,
with matrix representation $G_\theta$ (an $n\times n$ orthogonal matrix)
and $\theta=(\theta_1,\ldots,\theta_{\red{N}})$. 
The infinitesimal generators are
\begin{equation}\label{gj-def}
\g_j(Z) := \frac{\partial\ }{\partial\theta_j} G_\theta Z\bigg|_{\theta=0}\,,\quad
  j=1,\ldots,\red{N}\,.
\end{equation}
Since $G_\theta$ is orthogonal the action of $g_j$ on $Z$ is a skew-symmetric
matrix.
  Equivariance of (\ref{MJS}) then follows from the requirements
  \begin{equation}\label{MJS-equivariance}
  G_\theta{\bf M}={\bf M}G_\theta\,,\quad
  G_\theta{\bf J}={\bf J}G_\theta\,,\quad\mbox{and}\quad
  S(G_\theta Z)=S(Z)\,,\quad \forall\ G_\theta\in\T^{\red{N}}\,.
  \end{equation}

  The basic state, namely the solution that will be modulated,
  is taken to be a family of periodic $\red{N}-$phase wavetrains of the form
\begin{equation}\label{Zh-def}
Z(x,t) = \Zh(\bth,\bk,\bw)\,,\quad \bth = \bk x + \bw t  + \bth^{(0)}\,,
\end{equation}
\red{with the basic state
  (\ref{Zh-def}) $2\pi-$periodic
  in each component of $\bth$}, $\bth^{(0)}\in\mathbb{R}^N$ a constant,
and
\[
\bth = 
\begin{pmatrix}
\theta_1\\ \vdots \\
\theta_N
\end{pmatrix}\,, \quad
\bk = \begin{pmatrix}
k_1\\  \vdots \\
k_N
\end{pmatrix}\,, \quad
\bw = \begin{pmatrix}
\omega_1\\ \vdots \\
\omega_N
\end{pmatrix}\,.
\]
Substitution of $\Zh$ into (\ref{MJS}) admits the governing equation
for the $N-$phase wavetrain
\begin{equation}\label{BasicState}
\sum_{j=1}^N\big(\omega_j {\bf M}+k_j{\bf J}\big)\partial_{\theta_i}\Zh = \nabla S(\Zh)\,.
\end{equation}
In the absence of symmetry,
solutions of this problem may encounter small divisors.  
The advantage of the $\T^N-$symmetry of (\ref{MJS}) is that
multiphase wavetrains are smooth functions
with no small divisors: the wavetrain is a multiparameter family
of relative equilibria.
The relative equilibrium structure of the
basic state (\ref{Zh-def}) then gives
\begin{equation}\label{Zh-symmetry}
\Zh(\bth,\bw,\bk) = G_{\bth}\widehat{\bf z}(\bw,\bk)
\quad\mbox{with}\quad
\Zh_{\theta_j} = G_{\bth}\g_j(\widehat{\bf z})\,,\quad j=1,\ldots,N\,.
\end{equation}
All the dynamics is in the group action.  The
vector $\widehat{\bf z}(\bw,\bk)$
can be thought of as the reference point along the group orbit,
and satisfies
\[
\sum_{j=1}^N\big(\omega_j {\bf M}+k_j{\bf J}\big)\g_j(\widehat{\bf z}) = \nabla S(\widehat{\bf z})\,.
\]

\subsubsection{Linearization about a multiphase wavetrain}
\label{subsec-L-def}

Associated with (\ref{BasicState}) is the linear operator
\begin{equation}\label{L-def}
{\bf L} V = \D^2 S(\Zh) V - \sum_{j=1}^{\red{N}}\big(\omega _j{\bf M}+k_j{\bf J}\big)\partial_{\theta_j} V\,.
\end{equation}
This operator is formally self adjoint with respect to the inner product
\begin{equation} \label{ip-def}
\lth \cdot , \cdot \rth = \left(\frac{1}{2 \pi}\right)^{\red{N}}
\int_0^{2\pi}\cdots\int_0^{2\pi} \langle \cdot , \cdot \rangle\, \rd \theta _1 \cdots \rd \theta_{\red{N}} := \int_{\T^{\red{N}}}\langle\cdot,\cdot\rangle\, \rd\bth\,.
\end{equation}
Differentiation of (\ref{BasicState}) with respect to each
$\theta_i$ and each of the four parameters $k_i,\,\omega_i$ leads to the
equations
\begin{equation}\label{parameter-derivatives}
\begin{array}{rcl}
{\bf L}\Zh_{\theta_i} &=& 0\,, \\[2mm]
{\bf L}\Zh_{k_i} &=& {\bf J}\Zh_{\theta_i}\,,\\[2mm]
{\bf L}\Zh_{\omega_i} &=& {\bf M}\Zh_{\theta_i}\,,\quad i=1,\ldots,\red{N}\,.
\end{array}
\end{equation}
The first of these equations highlights the fact that the kernel
of ${\bf L}$ is at least
$N-$dimensional, and in this paper it is assumed no larger, so that
\begin{equation}\label{L-kernel}
  {\rm Ker}({\bf L}) = {\rm span}\left\{ \Zh_{\theta_1},\ldots,\Zh_{\theta_{\red{N}}} \right\}\,.
\end{equation}
The other equations in (\ref{parameter-derivatives})
will become significant when the Jordan chain theory in a moving frame is
developed.
The assumption (\ref{L-kernel}) along with the formal self-adjointness of
${\bf L}$ give the solvability conditions
for an expression $F$ to lie within the range of ${\bf L}$ as
\begin{equation}\label{L-solvability}
{\bf L}W = F \quad\Leftrightarrow\quad
\lth \Zh_{\theta_1},F\rth = \cdots = \lth \Zh_{\theta_{\red{N}}},F\rth = 0\,.
\end{equation}

\subsubsection{Multisymplectic Noether theory}
\label{subsec-claws}

In the Lagrangian setting, the
symmetry induces conservation laws via Noether theory.
Transforming to a multisymplectic formulation then induces
multisymplectic Noether theory which relates the structure operators
${\bf J}$ and ${\bf M}$ to the components of the induced
conservation laws.  Although these conservation laws may
have other physical significance they play the role of
conservation of wave action in the Whitham theory and so the
components will be called wave action and wave action flux.

There is a conservation law associated with each phase of the wavetrain,
and multisymplectic Noether theory implies the existence of functions $A_j,\,B_j$ satisfying
\begin{equation}
  {\bf M}\g_j(Z) = \nabla A_j(Z)\,, \quad {\bf J}\g_j(Z) = \nabla B_j(Z)\,,
  \quad j=1,\ldots,\red{N}\,,
\end{equation}
and so
\[
A_j(x,t) = \int_{\T^{\red{N}}}
\langle{\bf M}\g_j Z,Z\rangle\ d \bth\,,\quad B_j(x,t)
= \int_{\T^{\red{N}}}
\langle{\bf J}\g_j Z,Z\rangle\ d \bth\,,
\]
where $Z(x,t,\bth)$
is a function of $(x,t)$ and the phases
$\bth=(\theta_1,\ldots,\theta_{\red{N}})$ which are
here interpreted as ensemble parameters. 
Direct calculation verifies that the conservation laws are
\begin{equation}
\partial _t A_j +\partial _x B_j = 0\,, \quad j = 1,\ldots,\red{N}\,,
\end{equation}
whenever $Z$ satisfies (\ref{MJS}).   

The components of the conservation laws
can also be deduced directly from the averaged Lagrangian.
The Lagrangian (\ref{Lagr-mss}), 
evaluated on the $N-$phase wavetrain and averaged, is
\begin{equation}\label{avg-lagr}
\mathscr{L}(\bw,\bk) = \int_{\T^{\red{N}}}
\left[\sum_{j=1}^{\red{N}}\left[
 \fr\omega_j\langle \Zh, {\bf M}\Zh_{\theta_j} \rangle
+\fr k_j\langle \Zh , {\bf J} \Zh_{\theta_j} \rangle \right]
-S(\Zh) \right]\,\rd\bth\,.
\end{equation}
The wave action vector evaluated on the wavetrain is
\begin{equation}\label{Action-def}
{\bf A}(\bw,\bk) = \begin{pmatrix} \mathscr{A}_1 \\ \vdots\\
\mathscr{A}_{\red{N}}\end{pmatrix}
 := \D_\bw\mathscr{L}=\begin{pmatrix} \mathscr{L}_{\omega_1}\\
\vdots\\ \mathscr{L}_{\omega_{\red{N}}}\end{pmatrix} =\frac{1}{2}
\begin{pmatrix} \lth{\bf M}\Zh_{\theta_1},\Zh\rth \\ \vdots\\
 \lth{\bf M}\Zh_{\theta_{\red{N}}},\Zh\rth \end{pmatrix}\,,
\end{equation}
and the wave action flux vector is
\begin{equation}\label{Wave-action-def}
{\bf B}(\bw,\bk) = \begin{pmatrix} \mathscr{B}_1 \\  \vdots\\
  \mathscr{B}_{\red{N}}\end{pmatrix} := \D_\bk\mathscr{L}
  =\begin{pmatrix} \mathscr{L}_{k_1}\\ \vdots \\
\mathscr{L}_{k_{\red{N}}}\end{pmatrix} =\frac{1}{2}
\begin{pmatrix} \lth{\bf J}\Zh_{\theta_1},\Zh\rth \\ \vdots\\
 \lth{\bf J}\Zh_{\theta_{\red{N}}},\Zh\rth \end{pmatrix}\,.
\end{equation}
By definition, we have the following
\begin{equation}\notag
\begin{split}
&\hspace{2.1cm}\red{\D_\bw{\bf A} = \begin{pmatrix}
      \partial_{k_1}\mathscr{A}_1& \cdots & \partial_{k_{\red{N}}}\mathscr{A}_1\\
      \vdots & \ddots & \vdots \\
\partial_{k_1}\mathscr{A}_{\red{N}}&\cdots & \partial_{k_{\red{N}}}\mathscr{A}_{\red{N}}
\end{pmatrix}\,,}\\[4mm]
{\rm D}_{\bf k}{\bf A} =& \begin{pmatrix}
  \partial_{k_1}\mathscr{A}_1&\cdots&\partial_{k_{\red{N}}}\mathscr{A}_1\\
  \vdots & \ddots & \vdots \\
\partial_{k_1}\mathscr{A}_{\red{N}}&\cdots&\partial_{k_{\red{N}}}\mathscr{A}_{\red{N}}
\end{pmatrix} = {\rm D}_{\bm \omega}{\bf B}^T\,,\qand
{\rm D}_{\bf k}{\bf B} = \begin{pmatrix}
  \partial_{k_1}\mathscr{B}_1&\cdots&\partial_{k_{\red{N}}}\mathscr{B}_1\\
  \vdots & \ddots & \vdots \\
\partial_{k_1}\mathscr{B}_{\red{N}}&\cdots&\partial_{k_{\red{N}}}\mathscr{B}_{\red{N}}
\end{pmatrix}\,.
\end{split}
\end{equation}
The entries of these tensors are related to solutions via
\begin{subequations}
\begin{align}
\red{\partial _{\omega_j} \mathscr{A}_i} &= \red{\lth {\bf M} \Zh _{\theta _i} , \Zh _{\omega_j} \rth } \,, \label{Aw} \\
\partial _{k_j} \mathscr{A}_i &= \lth {\bf M} \Zh _{\theta _i} , \Zh _{k_j} \rth \,, \label{Ak} \\
\partial _{k_j} \mathscr{B}_i &= \lth {\bf J} \Zh _{\theta _i} , \Zh _{k_j} \rth \,, \label{Bk} \\
\partial _{k_j k_m} \mathscr{B} _i & = \lth {\bf J} \Zh_{\theta _i k_m} , \Zh _{k_j} \rth + \lth {\bf J} \Zh _{\theta _i} , \Zh _{k_j k_m} \rth\,,\quad i,j,m=1,\ldots,\red{N}\,.\label{Bkk}
\end{align}
\end{subequations}
The definition of the wave action and wave action flux in terms of
derivatives of the averaged Lagrangian induces symmetry of the Jacobians,
\begin{equation}\label{BDerivRelation}
\partial _{k_i} \mathscr{B} _j = \lth {\bf J} \Zh _{\theta _j} , \Zh _{k_i} \rth = \lth {\bf L} \Zh _{k _j} , \Zh _{k_i} \rth = \lth \Zh _{k _j} , {\bf L} \Zh _{k_i} \rth = \lth \Zh _{k _j} , {\bf J} \Zh _{\theta _i} \rth = \partial _{k_j} \mathscr{B}_i
\end{equation}
and
\begin{equation}\notag
  \partial _{k_j}\mathscr{A}_i = \lth {\bf M}\Zh_{\theta _i} , \Zh _{k_i} \rth = \lth \Zh _{\omega _i} , {\bf J} \Zh _{k _j}\rth = \partial _{\omega _i}\mathscr{B}_j
  \,,\quad i,j=1,\ldots,\red{N}\,.
\end{equation}
The key property in both (\ref{Action-def}) and (\ref{Wave-action-def}) is that
the left-hand side is in terms of the functions of $(\bw,\bk)$ only and the
right-hand side is expressed in terms of the properties of the Euler-Lagrange equation (\ref{MJS}), namely through the structure matrices ${\bf J}$ and ${\bf M}$. It is this connection that is the essence of multisymplectic Noether theory, and it feeds into the nonlinear modulation theory.

\subsection{Geometry of the averaged Lagrangian}
\label{sec-lagrangian-geometry}

Many of the properties needed in the modulation theory can be deduced
from the abstract mapping
\begin{equation}\label{L-avg}
  (\bw,\bk)\mapsto\mathscr{L}(\bw,\bk)\,,
\end{equation}
where $\mathscr{L}:\R^N\times\R^N\to\R$ is the averaged Lagrangian
(\ref{avg-lagr}) and is assumed to be a smooth function.

The wave action and wave action flux emerge from $\mathscr{L}$ via
\[
\frac{d\ }{ds} \mathscr{L}(\bw + s{\bf u} , \bk + s {\bf v})\Big|_{s=0} =
\langle{\bf A}(\bw,\bk),{\bf u}\rangle + \langle{\bf B}(\bw,\bk),{\bf v}\rangle\,,\quad \mbox{for any}\ {\bf u},{\bf v}\in\R^N\,,
\]
where $\langle\cdot,\cdot\rangle$ is an inner product on
$\R^N$.
The second derivative can be used to generate the linear operator
${\bf E}(c)$.  First set ${\bf u}=c{\bf v}$ in the above expression and look
at the derivative
\[
\begin{array}{rcl}
  \displaystyle
  \frac{d^2\ }{ds^2} \mathscr{L}(\bw + sc{\bf v} , \bk + s {\bf v})\Big|_{s=0}
  &=&\displaystyle
  \frac{d\ }{ds}
  \langle{\bf A}(\bw+sc{\bf v},\bk+s{\bf v}),c{\bf v}\rangle\Big|_{s=0}\\[3mm]
  &&\displaystyle \hspace{1.0cm}
  +  \frac{d\ }{ds}\langle{\bf B}(\bw+s{\bf v},\bk+s{\bf v}),{\bf v}\rangle\Big|_{s=0}\\[4mm]
  &=&\displaystyle
  \langle\D_\bw{\bf A}(\bw,\bk)c{\bf v},c{\bf v}\rangle
  +\langle\D_\bk{\bf A}(\bw,\bk){\bf v},c{\bf v}\rangle\\[3mm]  
  &&\displaystyle \hspace{1.0cm}
  +  \langle\D_\bw{\bf B}(\bw,\bk){\bf v},{\bf v}\rangle
   +  \langle\D_\bk{\bf B}(\bw,\bk){\bf v},{\bf v}\rangle\\[2mm]
&=&\displaystyle
   \left\langle \big(\D_\bw{\bf A}c^2 + (\D_\bk{\bf A}+\D_\bw{\bf B})c+\D_\bk{\bf B}\big){\bf v},{\bf v}\right\rangle  \\[3mm]
   &=&\displaystyle \langle{\bf E}(c){\bf v},{\bf v}\rangle\,,\quad\mbox{for
     any}\ {\bf v}\in\R^N\,.
\end{array}
\]
Hence
\begin{equation}\label{Ec-def-1}
  {\bf E}(c){\bf v} =
  \frac{d\ }{ds}\left[ c {\bf A}(\bw+cs{\bf v},\bk+s{\bf v})+
  {\bf B}(\bw+cs{\bf v},\bk+s{\bf v})\right]\bigg|_{s=0} 
 \,.
\end{equation}

The most remarkable result following from derivatives of $\mathscr{L}$
is the expression for $\kappa$, the coefficient of nonlinearity in
the emergent two-way Boussinesq equation (\ref{2-way-boussinesq-intro}).
Introduce the one parameter path in $\mathscr{L}(\bw,\bk)$,
\[
F(s) = \mathscr{L}(\bw + s c_g \be,\bk+s\be)\,,
\]
with $c_g$ here considered as fixed, and $\be\in{\rm Ker}({\bf E}(c_g))$.
  Then differentiating and using
(\ref{A-Lw-def}) and (\ref{B-Lk-def}) gives
\begin{equation}\label{F-Fp-Fpp}
\begin{array}{rcl}
F'(s) &=&\displaystyle \langle{\bf A}(\bw + s c_g \be,\bk+s\be),c_g\be\rangle +
\langle{\bf B}(\bw + s c_g \be,\bk+s\be),\be\rangle\\[2mm]
F''(s) &=&\displaystyle
\langle \D_\bw{\bf A}(\bw + s c_g \be,\bk+s\be)c_g\be,c_g\be\rangle
+\langle\D_\bk{\bf A}(\bw + s c_g \be,\bk+s\be)\be,c_g\be\rangle\\[2mm]
&&\hspace{1.0cm}\displaystyle
+ \langle\D_\bw{\bf B}(\bw + s c_g \be,\bk+s\be)c_g\be,\be\rangle
+ \langle{\bf B}_\bk(\bw + s c_g \be,\bk+s\be)\be,\be\rangle\,.\\[2mm]
\end{array}
\end{equation}
Evaluating $F''(0)$,
\[
F''(0) =
\langle c_g^2\D_\bw{\bf A}\be,\be\rangle
+\langle c_g\D_\bk{\bf A}\be,\be\rangle
+ \langle c_g\D_\bw{\bf B}\be,\be\rangle
+ \langle{\bf B}_\bk\be,\be\rangle =
\langle{\bf E}(c_g)\be,\be\rangle = 0 \,.
\]
However, it is the third derivative of $F(s)$ that is of most interest.
The formula for $F''(s)$ suggests that $F'''(s)$ is a derivative
of a path through the linear operator ${\bf E}(c_g)$, considered as a function
of $(\bw,\bk)$ with $c_g$ fixed.  Differentiating,
\begin{equation}\label{kappa-formula}
\begin{array}{rcl}
  F'''(0) &:=& \displaystyle
  \frac{d^3\ }{ds^3}\mathscr{L}(\bw + s c_g \be,\bk+s\be)\bigg|_{s=0}\\[4mm]
  &=& \left\langle\be,\big(\D_\bk^2{\bf B}+c_g(2\D_\bk\D_\bw{\bf B}+\D_\bk^2{\bf A})+c_g^2(2\D_\bk\D_\bw{\bf A}+\D_\bw^2{\bf B})+c_g^3\D_\bw^2{\bf A}\big)(\be,\be)
  \right\rangle\\[2mm]
  &:=& \kappa\,.
\end{array}
\end{equation}
At this point, this expression is just a formula, but the inner product
in the second row will emerge naturally in the modulation theory
in a solvability
condition, giving it relevance as the coefficient
of the nonlinear
term in the emergent modulation equation.

In a similar way, the
coefficient $\mu$ in (\ref{2-way-boussinesq-intro}) can also be
represented in terms of derivatives of $\mathscr{L}$ as in
\begin{equation}\label{alpha-def}
\mu=  \frac{d^2\ }{ds^2}\mathscr{L}(\bw + s\be,\bk)\bigg|_{s=0} +
\frac{d^2\ }{ds^2}\mathscr{L}(\bw+ sc_g{\bm\gamma},\bk+s{\bm\gamma})\bigg|_{s=0}\,.
\end{equation}
However, a more interesting characterization of $\mu$ is as a termination
condition for the Jordan chain $(\be,\bgam)$ in (\ref{zeta-def})-(\ref{gamma-def})
(see equations (\ref{Upsilon-eqn}) and (\ref{alpha-second-def}) below).

\section{Generic multiphase Whitham equations}
\setcounter{equation}{0}
\label{sec-generic-multiphase}

In this section a construction of the generic (distinct characteristics)
multiphase WMEs is sketched from the paper's geometric perspective.
  It serves as a touchstone for the
modifications needed for the non-generic (coalescing characteristics)
case, and the generic theory is needed to define $c_g$, the frame
speed at coalescence. 

Given the basic state $\Zh$ in (\ref{Zh-def}), \red{with $N$ phase variables $\theta_j$}, the generic WMEs
are obtained using the geometric
optics ansatz \cite{r18,r17b},
\begin{equation}\label{geom-optics-ansatz-1}
Z(x,t) = \Zh(\bth +\eps^{-1}\bph,\bw + \bO,\bk+\bq) + \eps W(\bth +\eps^{-1}\bph,X,T,\eps)\,
\end{equation}
with $X=\eps x$ and $T=\eps t$, and \red{the $N$-dimensional} vectors $\bph$,
$\bO$, and $\bq$  depending on $X,T$ and satisfying conservation
of waves $\bq_T=\bO_X$.
Expand all terms in a Taylor series, e.g. $W=W_1+\mathcal{O}(\eps)$,
substitute into (\ref{MJS}) and solve the equations at each order
of $\eps$.  At zeroth order the governing equations for the basic
wave $\Zh$
are recovered and at first order an equation for $W_1$ is
obtained
\[
   {\bf L}W_1 =  \red{\sum_{j=1}^N\bigg[\partial_T\Omega_j{\bf M}\Zh_{\omega_j}+\partial_T q_j {\bf M}\Zh_{k_j}+\partial_X\Omega_j{\bf J}\Zh_{\omega_j}+ \partial_Xq_j {\bf J}\Zh_{k_j}\bigg]\,.}
\]
Applying the solvability conditions (\ref{L-solvability}),
and using the connection between the resulting expressions
and the components of the conservation law (\ref{Aw})-(\ref{Bk}), \red{i.e}.
\[
\begin{array}{ccc}
\lth\Zh_{\theta_i},{\bf M}\Zh_{\omega_j}\rth = - \partial_{\omega_j}\mathscr{A}_i\,, &\lth\Zh_{\theta_i},{\bf J}\Zh_{k_j}\rth = - \partial_{k_j}\mathscr{A}_i\,,&\\[3mm]
\lth\Zh_{\theta_i},{\bf J}\Zh_{\omega_j}\rth = - \partial_{k\omega_j}\mathscr{B}_i\,, & \lth\Zh_{\theta_i},{\bf J}\Zh_{k_j}\rth = - \partial_{k_j}\mathscr{B}_i\,, & i,\,j=1,\ldots,\red{N}\,,
\end{array}
\]
then gives the generic WMEs,
\[
  0 = \red{\sum_{j=1}^N\bigg[\partial_T\Omega_j\partial_{\omega_j}\mathscr{A}_i+\partial_T q_j \partial_{k_1}\mathscr{A}_i +\partial_X\Omega_j\partial_{\omega_j}\mathscr{B}_i
  +\partial_X q_j\partial_{k_j}\mathscr{B}_i\bigg]\,,\quad i=1,\ldots,N\,.}
\]
Taking into account that $\Zh$ is a function of $\bk+\bq$ and
$\bw+\bO$, averaging over the phase eliminates the
$\eps^{-1}{\bm\phi}$ terms, and using the vector definition
of wave action (\ref{Action-def}) and wave action flux (\ref{Wave-action-def}),
these two equations are the
vector conservation equation
\begin{equation}\label{MP-WMEs-abstract-1}
  \partial_T{\bf A}(\bk+\bq,\bw+{\bm \Omega}) +
          \partial_X{\bf B}(\bk+\bq,\bw+{\bm \Omega}) =0\,,
\end{equation}
which, when combined with conservation of waves and the
symmetry condition
\begin{equation}\label{MP-WMEs-abstract-2}
\partial_T\bq=\partial_X{\bm \Omega} \qand \D_\bk{\bf A} = \big(\D_{\bw}{\bf B}\big)^T\,,
\end{equation}
give the generic WMEs in vector form. Further details of the above
derivation can be found in \cite{r18,r17b}.
A proof of validity of these multiphase WMEs, when the original equation is
  coupled NLS, covering both the cases of elliptic and hyperbolic
  characteristics, in the context of coupled nonlinear
  Schr\"odinger equations, is given in \textsc{Bridges et al.}~\cite{bks20}.

Consider the linearisation
of (\ref{MP-WMEs-abstract-1}) and (\ref{MP-WMEs-abstract-2}) at
$(\bw,\bk)$
\begin{equation}\label{wmes-linear}
  \D_\bw{\bf A}\,{\bm \Omega}_T+
  \D_\bk{\bf A}\,\bq_T+
  \D_\bw {\bf B}\,{\bm \Omega}_X+
  \D_\bk{\bf B}\,\bq_X={\bf 0}\qand
    {\bm \Omega}_X=\bq_T\,.
\end{equation}
    Characteristics about \emph{any} state $(\bw+\bO,\bk+\bq)$ can be obtained
    the same way, but here the main interest is in characteristics in the
    neighbourhood of the basic state. Differentiating the first equation and
    using the second results in a second order equation for $\bq$,
    \[
    \D_\bw{\bf A}\,\bq_{TT}+
  \big(\D_\bk{\bf A} + \D_\bw {\bf B}\big)\,\bq_{TX}+
  \D_\bk{\bf B}\,\bq_{XX}={\bf 0}\,.
  \]
With the normal mode ansatz
  \[
  ({\bm \Omega},\bq) = (\widehat{\bm \Omega},\widehat\bq)\re^{\ri\alpha(X+cT)}\,,
  \]
  the second-order equation results in a quadratic equation for the
  characteristics,
  \begin{equation}\label{wme-3-quadratic-a}
    {\bf E}(c)\widehat\bq:=
    \left[  \red{\D_\bw}{\bf A}\,c^2 + (\red{\D_\bk}{\bf A}+\red{\D_\bw}{\bf B})c + \red{\D_\bk}{\bf B}\right]
  \widehat\bq = {\bf 0}\,.
  \end{equation}
  It is a {\it Hermitian quadratic matrix polynomial}, and there is an extensive
  literature on the properties of these matrices
  (e.g. \textsc{Gohberg et al.}~\cite{glr80},
  \textsc{Tisseur \& Meerbergen}~\cite{tm01},
  \textsc{Mehrmann et al.}~\cite{mntx16} and references therein).
 
A key property that we will need is that a simple root, say $c_0$, has a
``sign characteristic''.  A necessary condition for two characteristics
to coalesce and transition
from hyperbolic to elliptic is that they have opposite sign
characteristic.  A study of the sign characteristic in the context of
the linearised multiphase WMEs is given in \cite{br19}, and the
basics of the theory needed here
is given below in Section \ref{sec-characteristics}.

\subsection{Pseudo-phases and affine symmetry}
\label{subsec-pseudo}

The basic states considered in \S\ref{subsec-symmetry}
  are $2\pi-$periodic in each phase (\ref{Zh-def}). There are
  also pseudo-phases.  In this section it is shown that pseudo-phases
  can be treated
  the same as phases associated with periodic motion, noting only that
  no averaging over pseudo-phases is required.
  A group-theoretic interpretation of
pseudo-phases is also given.

The concept of pseudo-phases was introduced by Whitham in his
first papers on modulation
(e.g.\ \textsc{Whitham}~\cite{whitham65,whitham67})
and discussed in more detail in \S14.6 in \cite{whitham-book}.
The associated pseudo frequencies and pseudo wavenumbers
play a role in mean flow. 

  Consider one of
  the most well-known examples of pseudo-phases; that is,
  mean flow in the water wave problem,
using the simplified model,
  \begin{equation}\label{ww-model-eqn}
  h_t + uh_x + hu_x=0 \qand u_t + uu_x + gh_x = \sigma h_{xxx}\,,
  \end{equation}
  where $h(x,t)$ is the fluid depth, $u(x,t)$ is the horizontal velocity
  field, and $\sigma$ is a parameter.
  Let $u=\phi_x$, and then the
  Lagrangian variational principle that generates (\ref{ww-model-eqn}) is
  \[
  \delta \int_{t_1}^{t_2}\int_{x_1}^{x_2} \mathcal{L}(\phi_t,\phi_x, h_t,h_x,h)\,
  \rd x\rd t=0\,,
  \]
  with
  \begin{equation}\label{Lagr-def}
  \mathcal{L}(\phi_t,\phi_x, h_t,h_x,h) = h\phi_t+
  \fr h\phi_x^2 + \fr g h^2 + \fr\sigma h_{x}^2 \,.
  \end{equation}
  The key property is that the Lagrangian is
  invariant if a constant is added to $\phi$.  This is an affine symmetry,
  the abstract group is the group of real numbers, and the action is
  \[
   s\cdot (h,\phi) = (h,\phi+s)\,,\quad \forall s\in\R\,.
  \]
This affine symmetry is reminiscent of ``cyclic variables'' in classical mechanics.

 When we introduce a pseudo-phase associated
 with this group, $s = \beta x+\gamma t$, the pseudo-phase does not
 appear explicitly in the Lagrangian and so averaging is not required.
 However, we will see that a modulation equation is still generated of
 the same mathematical form.  Let
 \begin{equation}\label{phi-pseudophase}
 \phi = \beta x + \gamma t\,,
 \end{equation}
 and substitute into (\ref{Lagr-def}),
  \begin{equation}\label{w-avg-lagr-0}
    \mathscr{L}(h,\beta) := \mathcal{L}(\gamma,\beta,0,0,h) =
    \gamma h + \fr h \beta^2 + \fr g h^2 \,.
  \end{equation}
  Although no averaging is required, Whitham modulation theory
  proceeds in the same way.  Suppose $\beta(X,T)$ and $\gamma(X,T)$
  are taken to depend on slow time and space variables.
  Then the WMEs associated with the single pseudo-phase are
  \begin{equation}\label{wmes-pseudophase}
  \beta_T = \gamma_X \qand \frac{\partial\ }{\partial T}\left(
  \frac{\partial\mathscr{L}}{\partial\gamma}\right) +
  \frac{\partial\ }{\partial X}\left(
  \frac{\partial\mathscr{L}}{\partial\beta}\right) =0\,,\quad
  \mbox{and}\quad \frac{\partial\mathscr{L}}{\partial h}=0\,.
  \end{equation}
  Differentiating $\mathscr{L}$ and substituting gives
  \[
  \beta_T = \gamma_X \qand \frac{\partial\ }{\partial T}\left(h
\right) +
  \frac{\partial\ }{\partial X}\left(\beta h\right)\,,\quad
  \mbox{and}\quad \gamma+\fr\beta^2 + gh=0\,.
  \]
  Substituting the third into the first, we get two equations
  \begin{equation}\label{swes}
  \beta_T + \beta\beta_X + gh_X =0\qand h_T + \beta h_X + h\beta_X=0\,.
  \end{equation}
  The upshot here is twofold, the introduction of a pseudo-phase does not
  require averaging, and the resulting modulation equations are mathematically
  the same as the case of modulation of periodic waves.

  Here we have the added outcome that Whitham theory applied
  to the pseudo-phase (\ref{phi-pseudophase})
  results in the classical shallow water equations (\ref{swes}),
  with $\beta(x,t)$ representing the horizontal velocity.

  If we add in a periodic phase as well, replacing (\ref{phi-pseudophase})
  with
  \[
  \phi(x,t) = \beta x + \gamma t + \Phi(\theta)\,,\quad
  \theta = kx+\omega t\,,
  \]
  where $\Phi$ is $2\pi-$periodic,
  and substitute into the Lagrangian, then averaging would be required but
  only over the periodic phase. Modulating $\omega(X,T)$ and $k(X,T)$,
  leads to the WMEs associated with the periodic phase
  \begin{equation}\label{wmes-phase}
  k_T = \omega_X \qand \frac{\partial\ }{\partial T}\left(
  \frac{\partial\mathscr{L}}{\partial\omega}\right) +
  \frac{\partial\ }{\partial X}\left(
  \frac{\partial\mathscr{L}}{\partial k}\right) =0\,,\quad
  \mbox{and}\quad \frac{\partial\mathscr{L}}{\partial E}=0\,,
  \end{equation}
  where $E$ is representative of the amplitude.
  Combining (\ref{wmes-pseudophase}) and (\ref{wmes-phase})
  results in coupled multiphase modulation equations with the
  same mathematical form as if both phases were
  periodic.  Indeed this is precisely what was done in
  \textsc{Whitham}~\cite{whitham67}.

  By eliminating the amplitudes in (\ref{wmes-pseudophase}) and (\ref{wmes-phase}) and relabelling $\gamma=\omega_1$, $\beta=k_1$, $\omega=\omega_2$, and
  $k= k_2$, the coupled multiphase modulation equations associated
  with pseudo-phase or phase take the same canonical form as (\ref{MP-WMEs-abstract-1}) and (\ref{MP-WMEs-abstract-2}). Hence the theory of this paper
  for multiphase Whitham modulation theory, including
  the theory for coalescing characteristics, applies to both periodic
  phases and pseudo-phases.  Henceforth the results will be stated for
  the periodic case and can easily be adjusted for pseudo-phases.
  
\section{Defining characteristics and coalescence}
\setcounter{equation}{0}
\label{sec-characteristics}

In this section the algebraic structure of the quadratic
Hermitian matrix pencil ${\bf E}(c)$ in (\ref{wme-3-quadratic-a})
is discussed.
Characteristics of the linearized WMEs (\ref{wmes-linear})
are the values of $c$ that are roots of the polynomial
\begin{equation}\label{Delta-def}
  \Delta(c) := {\rm det}[{\bf E}(c)]=0\,.
\end{equation}
When there are $N-$phases this polynomial has degree $2N$.
The linear algebra of quadratic Hermitian matrix pencils can be
found in \cite{glr80,mntx16,tm01} and references therein.
Here a theory for the sign characteristic of
simple roots and the theory of double non-semisimple roots
is required.

For definiteness we assume that all the characteristics are hyperbolic
and one pair transitions
from hyperbolic to elliptic at some parameter value.
It is not essential to the nonlinear theory for the 
other $\red{2(N-1)}$ characteristics to be hyperbolic at coalescence, although
if they are not hyperbolic then the basic state is already unstable.
      
      A characteristic is double when
      \begin{equation}\label{Delta-zero-eqn}
     \Delta(c_g)=\Delta'(c_g)=0 \qand \Delta''(c_g)\neq 0\,,
      \end{equation}
      where $\Delta(c)$ is defined in (\ref{Delta-def}).
The value of $c$ at the collision is denoted by $c=c_g$ in anticipation of
      the connection with the concept of group velocity. 

      The conditions (\ref{Delta-zero-eqn}) tell us that
      the algebraic multiplicity of $c_g$ is two.
      For Hermitian matrices the geometric multiplicity would also be two.
      However, Hermitian
     matrix pencils, in the indefinite case, can have non-trivial Jordan
     chains \cite{glr-book}.  This property also carries over to Hermitian
     quadratic matrix polynomials \cite{glr80}.  Here we are interested
     in the case where the
     geometric multiplicity of ${\bf E}(c_g)$ is one
     \begin{equation}\label{kernel-E}
       {\rm Ker}\big({\bf E}(c_g)\big) = {\rm span}\{\be\}\,.
     \end{equation}
     To establish a Jordan chain, first look at the condition
     $\Delta'(c_g)=0$ in terms of the properties of ${\bf E}(c)$,
\[
  \Delta'(c_g)
=\frac{d\ }{dc}{\rm det}[ {\bf E}(c)]\Big|_{c=c_g}\\[3mm]
  = {\rm Tr}\Big({\bf E}(c)^\#{\bf E}'(c)\Big)\Big|_{c=c_g}\,,
\]
where ${\bf E}(c)^\#$ is the adjugate \cite{mn88}.
Now use the fact that ${\bf E}(c_g)$ has rank one and the
nonzero eigenvalue is ${\rm Tr}({\bf E}(c_g))$,
  \[
    {\bf E}(c_g)^\# = \frac{{\rm Tr}({\bf E}(c_g))}{\|\be\|^2}\be\be^T\,.
    \]
    This formula can be verified by direct calculation (see also \cite{mn88}).
      Then
      \[
      \Delta'(c_g) =
             {\rm Tr}\Big({\bf E}(c)^\#{\bf E}'(c)\Big)\Big|_{c=c_g} =
             \frac{{\rm Tr}({\bf E}(c_g))}{\|\be\|^2}
             \langle\be,{\bf E}'(c_g)\be\rangle\,,
               \]
               and so with the assumption (\ref{kernel-E}),
               \begin{equation}\label{gen-eig-condition}
               \Delta'(c_g) =0 \quad\Longleftrightarrow\quad
               \langle\be,{\bf E}'(c_g)\be\rangle=0\,.
               \end{equation}
               Now look at this condition from the viewpoint of solvability, as that
               is how it will arise in the nonlinear modulation theory.
               In the case of algebraic multiplicity two and geometric
               multiplicity one, a Jordan chain for a quadratic Hermitian
               matrix polynomial has the form
     \begin{equation}\label{E-chain-two}
       {\bf E}(c_g)\be=0 \qand {\bf E}(c_g)\bgam = -{\bf E}'(c_g)\be\,,
     \end{equation}
     for some $\bgam\in\R^N$, if it exists \cite{glr80}.
     Since ${\bf E}(c_g)$ is
     Hermitian (in this case real and symmetric),
     the solvability condition is $\langle\be,{\bf E}'(c_g)\be\rangle=0$ confirming
     (\ref{gen-eig-condition}).  Writing out this condition,     
       \begin{equation}\label{E-solvability}
       0 = \langle\be,{\bf E}'(c_g)\be\rangle =
       \langle\be,\big(2c_g\D_{\bw}{\bf A} +(\D_{\bw}{\bf B}+\D_{\bk}{\bf A})\big)\be\rangle\,,
       \end{equation}
       gives a defining equation for $c_g$
\begin{equation}\label{cg-def}
  c_g = -\frac{1}{2}\frac{\langle\be,\big(\D_\bk{\bf A}+\D_\bw{\bf B}\big)\be\rangle}{\langle\be,\D_\bw{\bf A}\be\rangle}\,.
\end{equation}
Noting that $\D_\bw{\bf B} = (\D_\bk{\bf A})^T$, this formula simplifies to
\begin{equation}\label{cg-def-1}
  c_g = -\frac{\langle\be,\D_\bk{\bf A}\be\rangle}{\langle\be,\D_\bw{\bf A}\be\rangle}\,.
\end{equation}
The \red{notation} $c_g$ is used as the derivative with respect to $\bk$ over a
derivative with respect to $\bw$ is reminiscent of the classical definition of
group velocity.
               
Termination of the chain (\ref{E-chain-two})
at length two is assured if the following equation
\begin{equation}\label{Upsilon-eqn}
         {\bf E}(c_g)\Upsilon = -{\bf E}'(c_g)\bgam - \fr{\bf E}''(c_g)\be\,,
\end{equation}
         is not solvable; that is
         \begin{equation}\label{alpha-second-def}
         \mu := \langle\be,{\bf E}'(c_g)\bgam + \fr{\bf E}''(c_g)\be\rangle =\fr \langle\be,{\bf E}''(c_g)\be\rangle - \langle\bgam,{\bf E}(c_g)\bgam \rangle\neq 0\,,
         \end{equation}
         where (\ref{E-chain-two}) has been used.
         This expression is called $\mu$ as another remarkable result in
         the nonlinear theory is
         that this coefficient is precisely the $\mu$ that appears
         as the coefficient of $U_{TT}$ in the emergent two-way Boussinesq
         equation (\ref{2-way-boussinesq-intro}).  This connection will
         emerge in the nonlinear modulation theory.

         Further properties of this Jordan chain, and the Jordan chains
         associated with the linear operator ${\bf L}$ are discussed in
         more detail in \S\ref{sec-jordanchain}, after the nonlinear
         modulation theory is introduced.

\section{Nonlinear modulation at coalescence}
\label{sec-mod-ansatz}
\setcounter{equation}{0}

For the nonlinear modulation near coalescing characteristics,
the strategy is to introduce the ansatz (\ref{Z-cc-ansatz}),
substitute into the Euler-Lagrange equation (\ref{MJS}),
expand everything in powers of $\eps$, and set terms proportional to
each order in $\eps$ to zero.  The key step here is identifying the form
of the ansatz.  The role of frame speed is inspired by the one-phase case
in \cite{br17}, and the role of additional phase functions ${\bm\psi}$ and
${\bm\delta}$ is inspired by \cite{r18}.  These are included in the ansatz because they eliminate the need for
functions that would appear from homogeneous solutions at each order.
The proposed phase modulation is
\begin{equation}\label{finer-ansatz-Phi}
  \Phi = \bph+\eps \bps +\eps^2{\bm\delta}\,.
\end{equation}
Then with
\begin{equation}\label{finer-ansatz-omega-k}
  \bO := \bph_T\qand \bq := \bph_X\quad\Rightarrow\quad
  \bq_T - \bO_X =0\,,
\end{equation}
and
\begin{equation}\label{xt-cg-def}
X = \eps (x + c_g t)\,,\quad T=\eps^2t\,,
\end{equation}
the complete proposed ansatz (\ref{Z-cc-ansatz}) is
\begin{equation}\label{Ansatz}
\begin{array}{rcl}
Z(x,t) &=& \Zh\big(\bth+{\eps}\bph+{\eps^2}\bps+{\eps^3}{\bm \delta},\bk+{\eps^2}\bq+{\eps^3}\bps_X+{\eps^4}{\bm \delta}_X,\\[2mm]
&&\hspace{1.0cm} \bw+{\eps^2}c_g\bq+{\eps^3}(\bO+c_g\bps_X)+{\eps^4}(\bps_T+c_g{\bm \delta_X})+{\eps^5}{\bm \delta}_T\big)\\[2mm]
&&\hspace{2.0cm} +{\eps^3}W(\bth,X,T;{\eps})\,.
\end{array}
\end{equation}
where $\bth$, $\bph$, $\bps$, ${\bm\delta}$, $\bq$, and $\bO$
are all functions of $X$ and $T$ defined in
(\ref{xt-cg-def}), and $c_g$ is defined in (\ref{cg-def}).
For ease, expand $W$ in an asymptotic series, 
\[
W(\bth,X,T,{\eps}) = W_3(\bth,X,T)
+ {\eps} W_4(\bth,X,T)+
{\eps^2}W_5(\bth,X,T)+\ldots\,.
\]
The remainder $W$ could be defined as
$W(\bth+{\eps}\bph+{\eps^2}\bps+{\eps^3}{\bm \delta},X,T,\eps)$, to synchronise
with the form of the modulation of the basic state, but is equivalent to
the above formulation: expansion of $W$ in a Taylor series in $\eps$
just changes the form of $W_j$ at each order, but the overall expansion
gives equivalent results.

Although the ansatz (\ref{Ansatz}) is new, the expansion
and substitution strategy is similar to our
previous papers on multiphase modulation \cite{rb16a,r17b,r18}
and the single phase coalescing characteristics \cite{br17}
and so only the key new points are highlighted.
For example, at $\eps^0$ order the governing equation for $\Zh$ in (\ref{BasicState})
is recovered.  At $\eps^1$ and $\eps^2$ order the generic $2-$term
Jordan chain in (\ref{parameter-derivatives}) is recovered as in the
preceding works.

At third order in $\eps$, after simplification, the system is
\begin{equation}\label{W3-eqn}
  {\bf L}W_3 = \sum_{j=1}^{\red{N}} \partial_X q_j {\bf K}\left(\Zh_{k_j}+c_g\Zh_{\omega_j}\right)\,,
\end{equation}
where
\begin{equation}\label{K-def}
{\bf K} := {\bf J}+c_g {\bf M}\,,
\end{equation}
and $c_g$ is defined in (\ref{cg-def}). 
Applying the solvability condition (\ref{L-solvability})
to (\ref{W3-eqn}) gives 
\[
\begin{array}{rcl}
 \displaystyle \sum_{j=1}^{\red{N}}\partial_X q_j\blth\Zh_{\theta_i}, ({\bf J}+c_g{\bf M})\left(\Zh_{k_j}+c_g\Zh_{\omega_j}\right)\brth &=& 0\,, \quad \red{i = 1,\ldots,N}
\end{array}
\]
or, after using the conversions from the structure operators
${\bf J}$, ${\bf M}$ to the functionals $\mathscr{A}_j,\mathscr{B}_j$ in
(\ref{Ak})-(\ref{Bk}), the solvability condition can be written in the
illuminating vector form
\begin{equation}\label{E-solve-q}
\left[\D_\bk{\bf B}+c_g(\D_\bw{\bf B}+\D_\bk{\bf A})+c_g^2\D_\bw{\bf A}\right]\bq_X = {\bf 0}\,.
\end{equation}
Hence for solvability of (\ref{W3-eqn})
it is required that $\bq_X$ is in the kernel
of ${\bf E}(c_g)$,
\[
\bq_X = U_X\be\quad\Rightarrow\quad \bq = U(X,T)\be + {\bf a}(T)\,,
\]
for some scalar-valued function $U(X,T)$.  It can be confirmed
\emph{a posteriori} that ${\bf a}(T)$ does not contribute to the leading
order result and can be neglected.  Hence
\begin{equation}\label{q-U-def}
\bq = U(X,T)\be\,.
\end{equation}
It is this scalar-valued function $U(X,T)$ that will ultimately
be found to be governed by the
two-way Boussinesq equation (\ref{2-way-boussinesq-intro}).

With the solvability condition satisfied, and the expression for
${\bf q}$ in (\ref{q-U-def}), the complete solution at third order is
\begin{equation}\label{W3-soln}
W_3 = U_X{\bf v}_3\,, \quad\mbox{with}\quad
{\bf L}{\bf v}_3 = {\bf K}{\bf v}_2\,.
\end{equation}
An arbitrary amount of homogeneous solution can be added to $W_3$ but it
is already incorporated into the functions ${\bm\delta}$ and ${\bm\psi}$ in the ansatz.
The equation ${\bf Lv}_3={\bf Kv}_2$ foreshadows a Jordan chain theory.
The beginnings of the chain are in (\ref{parameter-derivatives})
which can be re-written as ${\bf Lv}_1=0$ and ${\bf Lv}_2={\bf Kv}_1$.
This Jordan chain theory is developed in Section \ref{sec-jordanchain}.

\subsection{Fourth order}
\label{subsec-fourthorder}

After simplification, the equation at fourth order is
\begin{equation}\label{LW4-def}
\begin{array}{rcl}
  \displaystyle{\bf L}\left( W_4 - U_{X}\sum_{i=1}^{\red{N}} \phi_i({\bf v}_3)_{\theta_i}\right)
  &=& \displaystyle
  U_{XX}{\bf K}{\bf v}_3  +\sum_{j=1}^{\red{N}}(\psi_j)_{XX}{\bf K}\left(\Zh_{k_j}+c_g\Zh_{\omega_j}\right)\\[2mm]
  &&\displaystyle\quad
  + U_T\sum_{j=1}^{\red{N}}\zeta_j\big({\bf J}\Zh_{\omega_j}
  +{\bf M}\Zh_{k_j}+2c_g{\bf M}\Zh_{\omega_j}\big)\,.
\end{array}
\end{equation}
The first inhomogeneous term feeds into the Jordan chain argument as it is
of the form ${\bf L}{\bf v}_4={\bf K}{\bf v}_3$, for some ${\bf v}_4$.
For the other two inhomogeneous terms, apply the solvability conditions
(\ref{L-solvability}), and use the identities (\ref{Ak})-(\ref{Bk}), to obtain
\begin{equation}\label{solv-order-4}  
{\bf E}(c_g)\bps_{XX}+\underbrace{\left[(\D_\bk{\bf A}+\D_\bw{\bf B})+2c_g\D_\bw{\bf A}\right]}_{{\bf E}'(c_g)}\be U_T  = {\bf 0}\,.
\end{equation}
This equation is of the form (\ref{E-chain-two}); that is, the
Jordan chain associated with ${\bf E}(c)$.  The theory of this
Jordan chain is developed below in
\S\ref{subsec-cg-jordanchain}.  Here, it is sufficient to
use the argument presented in (\ref{E-chain-two}) and (\ref{cg-def})
for the chain $(\be,\bgam)$ of ${\bf E}(c_g)$.  Applying that
theory gives
\begin{equation}\label{psi-XX-def}
\bps_{XX} =\bgam U_T\quad \mbox{(mod Ker(${\bf E}(c_g)$))}\,,
\end{equation}
where ``mod'' signifies that an arbitrary amount of homogeneous solution
can be included.  This homogeneous solution can
be neglected as it does not enter at fifth order.
Thus the solution at fourth order is
\begin{equation}\label{W4-soln}
  W_4 = U_{XX}{\bf v}_4+U_T\Xi+U_X \sum_{j=1}^{\red{N}}\phi_j(\xi_5)_{\theta_j}
  \quad \mbox{(mod Ker(${\bf L}$))}\,,
\end{equation}
with
\begin{equation}\label{Xi-def-eqn}
{\bf L}\Xi = \sum_{j=1}^{\red{N}}\left[\zeta_j\big({\bf J}\Zh_{\omega_j}+{\bf M}\Zh_{k_j}+2c_g{\bf M}\Zh_{\omega_j}\big)+\gamma_j{\bf K}(\Zh_{k_j}+c_g\Zh_{\omega_j})\right]\,.
\end{equation}
Fortunately this equation does not need to be solved explicitly.  It feeds
into the fifth order solution, and ultimately generates formulae for the
coefficients, but these formulae will be obtained without an explicit
expression for $\Xi$.

\subsection{Fifth order}
\label{subsec-fifthorder}

At fifth order, after combining terms and simplifying, the equations are
\begin{equation}\label{fifthorder}
\begin{array}{rcl}
  {\bf L}\widetilde{W_5} &=& \displaystyle
  U_{XXX}{\bf K}{\bf v}_4+\sum_{i=1}^{\red{N}}\bigg[(\Omega_i)_T{\bf M}\Zh_{\omega_i}+(\delta_i)_X{\bf K}(\Zh_{k_i}+c_g\Zh_{\omega_i})\bigg]\\[4mm]
  &&\displaystyle\quad
+U_{XT}({\bf J}\Xi+{\bf M}{\bf v}_3)
  +\sum_{i=1}^{\red{N}} (\psi_i)_{XT}( {\bf J}\Zh_{\omega_i}+{\bf M}\Zh_{k_i}+
  2c_g{\bf M}\Zh_{\omega_i})\\[2mm]
&&\displaystyle\quad + UU_X\sum_{i=1}^{\red{N}}\bigg[{\bf K}({\bf v}_3)_{\theta_i}-\D^3S(\Zh)({\bf v}_3,\Zh_{k_i}+c_g\Zh_{\omega_i})\\[2mm]
    &&\displaystyle
    \hspace{3.5cm}\displaystyle +\sum_{j=1}^{\red{N}}{\bf K}(\Zh_{k_ik_j}+c_g\Zh_{\omega_ik_j}+c_g\Zh_{k_i\omega_j}+c_g^2\Zh_{\omega_i\omega_j})\bigg]\,.
\end{array}
\end{equation}
The tilde above $W_5$ term indicates that
the preimage of all terms lying in the range of ${\bf L}$ from the
right-hand side have been absorbed (e.g.\ terms that would vanish
identically under the solvability conditions).
These terms come into play only at higher order.

It is the solvability condition for this
fifth order equation that will deliver the modulation equation for $U(X,T)$.
However, solvability is a multistage process.  There are \red{$N$} solvability
conditions associated with the operator ${\bf L}$, leading to a vector-valued
equation.  A secondary solvability condition, associated with the
operator ${\bf E}(c_g)$, will reduce vector equation to the scalar
two-way Boussinesq equation.

First establish the vector solvability condition.
Apply the ${\bf L}-$solvability (\ref{L-solvability}) condition
to the right-hand side of (\ref{fifthorder})
term by term.
Solvability of the $U_{XXX}$ term generates the vector
\begin{equation}\label{DispTerm}
\begin{pmatrix}
\lth \Zh_{\theta_1},{\bf K}{\bf v}_4\rth\\
\red{\vdots}\\
\red{\lth \Zh_{\theta_N}},{\bf K}{\bf v}_4\rth
\end{pmatrix}U_{XXX} :=-{\bf T}U_{XXX}\,.
\end{equation}
We will see that this vector is nonzero since, by hypothesis, the Jordan chain
$({\bf v}_1,\ldots,{\bf v}_4)$ has length four.  This is discussed below
in \S\ref{sec-jordanchain}.
Solvability of the $(\Omega_i)_T$ terms leads to the matrix term
\begin{equation}\label{TimeTerm1}
\begin{pmatrix}
\lth \Zh_{\theta_1},{\bf M}\Zh_{\omega_1}\rth&\red{\cdots}&\red{\lth \Zh_{\theta_1},{\bf M}\Zh_{\omega_N}\rth}\\
\vdots&\ddots&\vdots\\
\lth \Zh_{\theta_{\red{N}}},{\bf M}\Zh_{\omega_1}\rth&\red{\cdots}&\lth \Zh_{\theta_{\red{N}}},{\bf M}\Zh_{\omega_{\red{N}}}\rth
\end{pmatrix}\bO_T \equiv -\D_\bw {\bf A}\,\bO_T\,.
\end{equation}
The terms containing $\delta_i$ give
\begin{equation}\label{AlphaTerm}
\begin{pmatrix}
\lth \Zh_{\theta_1},{\bf K}(\Zh_{k_1}+c_g\Zh_{\omega_1})\rth&\red{\cdots}&\lth \Zh_{\theta_1},{\bf K}(\Zh_{k_{\red{N}}}+c_g\Zh_{\omega_{\red{N}}})\rth\\
\red{\vdots}&\red{\ddots}&\red{\vdots}\\
\lth \Zh_{\theta_{\red{N}}},{\bf K}(\Zh_{k_1}+c_g\Zh_{\omega_1})\rth&\red{\cdots}&\lth \Zh_{\theta_{\red{N}}},{\bf K}(\Zh_{k_{\red{N}}}+c_g\Zh_{\omega_{\red{N}}})\rth
\end{pmatrix}{\bm \delta}_{XX} = -{\bf E}(c_g){\bm \delta}_{XX}\,.
\end{equation}
The terms involving $(\psi_i)_{XT}$ are similar to those seen at fourth order,
and generate
\begin{equation}\label{Timeterm2}
\begin{split}
\begin{pmatrix}
\lth \Zh_{\theta_1},{\bf J}\Zh_{\omega_1}+{\bf M}\Zh_{k_1}+2c_g{\bf M}\Zh_{\omega_
1}\rth&\red{\cdots}&\lth \Zh_{\theta_1},{\bf J}\Zh_{\omega_{\red{N}}}+{\bf M}\Zh_{k_{\red{N}}}+2c_g{\bf M}\Zh_{
\omega_{\red{N}}}\rth\\
\red{\vdots}&\red{\ddots}&\red{\vdots}\\
\lth \Zh_{\theta_{\red{N}}},{\bf J}\Zh_{\omega_1}+{\bf M}\Zh_{k_1}+2c_g{\bf M}\Zh_{\omega_
1}\rth&\red{\cdots}&\lth \Zh_{\theta_{\red{N}}},{\bf J}\Zh_{\omega_{\red{N}}}+{\bf M}\Zh_{k_{\red{N}}}+2c_g{\bf M}\Zh_{
\omega_{\red{N}}}\rth
\end{pmatrix}
\bps_{XT}\\
= -\bigg[(\D_\bk{\bf A}+\D_\bw{\bf B})+2c_g\D_\bw{\bf A}\bigg]\bps_{XT}\\
= -{\bf E}'(c_g)\bps_{XT}\,.
\end{split}
\end{equation}
The coefficient of the nonlinear term $UU_X$ simplifies to
\begin{equation}\label{QuadTerm}
  \begin{split}-(\D_\bk^2{\bf B}+c_g(2\D_\bk\D_\bw{\bf B}+\D_\bk^2{\bf A})+c_g^2(2\D_\bk\D_\bw{\bf A}+\D_\bw^2{\bf B})+c_g^3\D_\bw^2{\bf A})(\be,\be)UU_X\\
    := -{\bf H}(\be,\be) UU_X\,.
    \end{split}
\end{equation}
When $c_g=0$ the vector function ${\bf H}(\be,\be)$ reduces to
${\bf H}(\be,\be) = \D^2_{\bk}{\bf B}(\be,\be)$ which is the form found
in reduction of multiphase modulation to the KdV equation
in \cite{rb16a,r18b}.

Collecting these terms gives the vector form of the solvability condition
for (\ref{fifthorder})
\begin{equation}\label{vector-fifth-order-equation}
  {\bf E}(c_g){\bm \delta}_{XX} + \D_\bw {\bf A}\Omega_T
+ {\bf E}'(c_g){\bm\psi}_{XT}
+ {\bf T}U_{XXX} +{\bf H}(\be,\be) UU_X = {\bf 0}\,.
\end{equation}
This equation is interesting in itself, but it is not closed due
to the presence of the ${\bm\delta}_{XX}$ term
and the ${\bm\psi}_{XT}$ term.  However, the ${\bm\delta}_{XX}$ term is
acted on by $\red{{\bf E}(c_g)}$ and so this term vanishes when the equation
is projected onto the kernel of ${\bf E}(c_g)$.  Therefore
split $\R^{\red{N}}$ as
\[
\R^{\red{N}} = {\rm span}\{\be\}\oplus \R^{\red{N-1}}\,.
\]
The projection of (\ref{vector-fifth-order-equation}) onto the
complement of ${\rm Ker}({\bf E}(c_g))$ still contains the
${\bm\delta}_{XX}$ term but this part carries over to higher order.

With this splitting in mind, act on (\ref{vector-fifth-order-equation}) with $\be^T$,
\begin{equation}\label{scalar-fifth-order-equation}
  \be^T{\bf E}(c_g){\bm \delta}_{XX} + \be^T\D_\bw {\bf A}\bO_T
+ \be^T{\bf E}'(c_g){\bm\psi}_{XT}
+ \be^T{\bf T}U_{XXX} +\be^T{\bf H}(\be,\be) UU_X = 0\,.
\end{equation}
Defining
\[
\kappa = \be^T{\bf H}(\be,\be)\qand \mathscr{K}=
\be^T{\bf T} = \lth{\bf Kv}_1,{\bf v}_4\rth\,,
\]
and noting that the coefficient of the ${\bm \delta}_{XX}$ term now vanishes 
as ${\bf E}(c_g)$ is symmetric,
(\ref{scalar-fifth-order-equation}) simplifies the vector equation to
\begin{equation}\label{pre-2way-Boussinesq-equation}
 \be^T\D_\bw {\bf A}\bO_T
+ \be^T{\bf E}'(c_g){\bm\psi}_{XT}
 + \kappa UU_X + \mathscr{K}U_{XXX} = 0\,.
\end{equation}
This equation is closed by first differentiating with respect to $X$,
\[
\be^T\D_\bw {\bf A}\bO_{XT}
+ \be^T{\bf E}'(c_g){\bm\psi}_{XXT}
 + \kappa (UU_X)_X + \mathscr{K}U_{XXXX} = 0\,,
 \]
 and applying conservation of waves and the
 ${\bm\psi}$-$U$ equation (\ref{psi-XX-def}),
\[
\bO_{XT} = \bq_{TT} = \be U_{TT} \qand {\bm\psi}_{XXT}={\bm\gamma} U_{TT}\,.
\]
Hence the final form of the
two-way Boussinesq equation is
\begin{equation}\label{U-2way-Boussinesq}
  \mu U_{TT}  + \kappa\, (UU_X)_X + \mathscr{K} U_{XXXX}=0\,,
\end{equation}
with
\begin{equation}\label{alpha-redef}
\mu =\be^T\D_\bw {\bf A}\be 
+ \be^T\left[ (\D_{\bk}{\bf A}+\D_{\bw}{\bf B}-2c_g\D_{\bw}{\bf A}\right] {\bm\gamma}\,.
\end{equation}
Another way to write this is to use ${\bf E}(c_g)$,
\begin{equation}\label{alpha-redef-1}
\mu =\fr\be^T {\bf E}''(c_g)\be 
+ \be^T {\bf E}'(c_g) {\bm\gamma}\,;
\end{equation}
emphasising that $\mu\neq0$ is the termination condition for the $(\be,\bgam)$
Jordan chain.

Comparing $\zeta^T{\bf H}(\be,\be)$ with (\ref{kappa-formula}) shows that
\begin{equation}\label{kappa-redef}
\kappa = \be^T{\bf H}(\be,\be) = \frac{d^3\ }{ds^3}\mathscr{L}(\bw + s c_g \be,\bk+s\be)\bigg|_{s=0}\,.
\end{equation}
The emergent two-way Boussinesq equation is non-degenerate when $\mu$,
$\kappa$ and $\mathscr{K}$ are nonzero.  The coefficient $\mu$ is nonzero
when the Jordan chain for ${\bf E}(c_g)$ in (\ref{E-chain-two}) terminates
at two.  The coefficient $\kappa$ is assumed to be nonzero.  If it is
zero, then it is expected that re-modulation will lead to a cubic
nonlinearity \cite{ehs17,rb18,r18b}.
The
coefficient of dispersion is nonzero if the Jordan chain 
in (\ref{four-chain}) terminates at four.  If $\mathscr{K}$
vanishes, then a longer Jordan chain will emerge.  Re-modulation in
this case is expected to lead to higher order dispersive terms emerging
(e.g.\ sixth-order dispersion, as in \cite{sh17,r17c}).

The above result does not provide any information about
convergence of the ansatz (\ref{Ansatz}) as a Taylor series in $\eps$.
However, the \emph{asymptotic validity} of this ansatz is confirmed by
the above results; that is,
the ansatz (\ref{Ansatz}) satisfies the governing equations
exactly up to $\mathcal{O}(\eps^5)$,
\[
\Big\| {\bf M}Z_t+{\bf J}Z_x = \nabla S(Z)\Big\| = \mathcal{O}(\eps^6)\quad\mbox{as}\ \eps\to0\,.
\]
For generic multiphase WMT, a rigorous proof of validity has been given
for CNLS \cite{bks20}, but a rigorous proof of validity in the case of
coalescing characteristics is an open problem.

To summarize, the starting point is a PDE generated by a Lagrangian with
a multiphase basic state.  It is assumed that, at some parameter value,
a pair of
coalescing characteristics arises in the linearized Whitham equations.
These coalescing characteristics generate several Jordan chains.
A modulation ansatz of the form (\ref{Ansatz}) then leads to
a scalar two-way Boussinesq equation (\ref{U-2way-Boussinesq})
with coefficients $\mu$, $\kappa$, and $\mathscr{K}$ all determined
from abstract properties of the averaged Lagrangian.
The fundamental idea is
that the original PDE is reduced to a simpler PDE that can be analyzed in some
detail.  Some of the solutions of this reduced two-way Boussinesq equation are
anticipated in \S\ref{sec-boussinesq-eqn}.

\section{Coalescing characteristics and Jordan chains}
\label{sec-jordanchain}
\setcounter{equation}{0}

Jordan chains play an important part throughout the
steps of the derivation of the nonlinear modulation equations.
In this section, some of the properties
of these Jordan chains are examined in more detail.

There are two key linear operators: ${\bf L}$ and
${\bf E}(c)$.  The operator ${\bf L}$, associated with the linearization of
the Euler-Lagrange equation (\ref{MJS}), generates a Jordan chain theory
that starts with
\begin{equation}\label{xi1-xi2-def}
  {\bf L}\xi_j = 0\quad \mbox{with}\quad
  \xi_j := \frac{\partial\Zh}{\partial\theta_j}\,,\ j=1,\ldots,\red{N}\,,
\end{equation}
and
  \begin{equation}\label{J-chains}
  \begin{array}{rcclcl}
    {\bf L}\xi_1 &=& 0\,, & {\bf L} \xi_{\red{N}+1} &=& {\bf J}\xi_1\\[2mm]
    \vdots   &&  &&  \vdots & \\
             {\bf L}\xi_{\red{N}} &=& 0\,, & {\bf L} \xi_{2\red{N}}
             &=& {\bf J}\xi_{\red{N}}\,,
                 \end{array}
  \end{equation}
    which follow from (\ref{parameter-derivatives}) with
    \[
    \xi_{\red{N}+j} :=\frac{\partial\Zh}{\partial k_j}\,,\quad j=1,\ldots,\red{N}\,.
    \]
With the assumption that ${\rm Ker}({\bf L})={\rm span}\{\xi_1,\ldots,\xi_{\red{N}}\}$,
there are $\red{N}$ Jordan blocks each of dimension two.
    
    The operator ${\bf E}(c_g)$ generates another Jordan chain which can
    be discussed independently
    of the ${\bf L}-$chains, but feeds into solvability of the
    ${\bf L}$ chains, and it starts with
    \begin{equation}\label{E-kernel}
      {\bf E}(c)\be = 0 \,.
    \end{equation}
    Generically, ${\bf E}(c)$ has a single Jordan block of dimension one.
    \red{At isolated values, for example at $c=c_g$, the Jordan block increases
    to dimension two.  Even though ${\bf E}(c)$ is Hermitian for any real
    $c$, it still generates
    a non-trivial Jordan chain due to the fact that $c$ appears nonlinearly.}
    
    The theory needed to extend these two Jordan chains is well
    established in the literature.  The above ${\bf L}-$chains are
    ${\bf J}-$symplectic Jordan chains and this theory goes back
    to \textsc{Williamson}~\cite{williamson}, and the theory
    of Jordan chains for quadratic Hermitian matrix pencils is developed
    in \cite{glr80}.

    However, things get complicated when we realise
    that the linear operator ${\bf L}$ has both ${\bf J}-$Jordan
    chains and ${\bf M}-$Jordan chains.  From
    (\ref{parameter-derivatives}) it follows that there exist
    ${\bf M}-$Jordan chains of the form
  \begin{equation}\label{M-chains}
  \begin{array}{rcclcl}
    {\bf L}\xi_1 &=& 0\,, & {\bf L} \eta_{\red{N}+1} &=& {\bf M}\xi_1\\[2mm]
    \vdots   &&  &&  \vdots & \\
                 {\bf L}\xi_{\red{N}} &=& 0\,, & {\bf L} \eta_{2\red{N}} &=& {\bf M}\xi_{\red{N}}\,,
                 \end{array}
  \end{equation}
    which follow from (\ref{parameter-derivatives}) with
    \[
    \eta_{\red{N}+j} :=\frac{\partial\Zh}{\partial k_j}\,,\quad j=1,\ldots,\red{N}\,.
    \]    
    The ${\bf J}-$chains (\ref{J-chains}) have length greater than two
    if at least one of the following
    \[
      {\bf L}\chi_j ={\bf J}\xi_{\red{N}+j}\,,\quad j=1,\ldots,\red{N}\,,
    \]
      is solvable, and termination at two is associated with
      non solvability of all $\red{N}$ of these equations.
      Similarly, the ${\bf M}-$chains have length greater than two if at least
      one of the following
    \[
      {\bf L}\chi_j ={\bf M}\eta_{\red{N}+j}\,,\quad j=1,\ldots,\red{N}\,,
        \]
    is solvable, and termination at two is associated with
    non solvability.
    The $\red{N}$ chains in (\ref{J-chains}) and (\ref{M-chains})
      can also be mixed, by
      taking the first elements to be linear combinations of
      $\xi_1,\ldots,\xi_{\red{N}}$, and this turns out to be useful for
      the modulation theory.

      Combining all the possibilities for both ${\bf J}-$chains
      and ${\bf M}-$chains, the most
      general extension of the Jordan chains is
that there exists a vector $\Xi$ satisfying
\begin{equation}\label{L.1}
  {\bf L}\Xi
  = \sum_{j=1}^{\red{N}} \left( a_j{\bf M}\eta_{\red{N}+j}  + b_j{\bf M}\xi_{\red{N}+j}  
   +c_j{\bf J}\eta_{\red{N}+j}  + d_j{\bf J}\xi_{\red{N}+j} \right)\,.
\end{equation}
No theory exists for Jordan chains of this type.  The closest
approximation is the Jordan chain theory for multiparameter eigenvalue
problems (e.g.\ \textsc{Binding \& Volkmer}~\cite{bv96} and
its citation trail),
but that does not apply here either.  We will be able to develop
a satisfactory theory for multi-dimensional Jordan chains of this type
to cover the cases needed in the modulation theory, but a complete
and general theory for multi-dimensional Jordan chains of this type
is outside the scope of this paper.

A solution $\Xi$ of (\ref{L.1}) exists if this equation is solvable,
and it is solvable
if and only if the $4\red{N}$ constants $a_j,b_j,c_j,d_j$ for $j=1,\ldots,\red{N}$,
satisfy
\[
\left\langle \xi_\ell,
\sum_{j=1}^{\red{N}} \left( a_j{\bf M}\eta_{\red{N}+j}  + b_j{\bf M}\xi_{\red{N}+j}  
+c_j{\bf J}\eta_{\red{N}+j}  + d_j{\bf J}\xi_{\red{N}+j} \right)\right\rangle=0\,,
\quad \ell=1,\ldots,\red{N}\,.
\]
These equations can be simplified
by using the identities in \S\ref{subsec-claws}, giving
  \begin{equation}\label{E-solvability-2}
  [\D_{\bw}{\bf A}] {\bf a} +
  [\D_{\bk}{\bf A}] {\bf b} +
  [\D_{\bw}{\bf B}] {\bf c} +
  [\D_{\bk}{\bf B}] {\bf d}   = {\bf 0}\,,
  \end{equation}
  where
  \[
  {\bf a}:=\begin{pmatrix} a_1\\ \vdots\\ a_{\red{N}}\end{pmatrix}\,,\quad
  {\bf b}:=\begin{pmatrix} b_1\\ \vdots\\ b_{\red{N}}\end{pmatrix}\,,\quad
  {\bf c}:=\begin{pmatrix} c_1\\ \vdots\\ c_{\red{N}}\end{pmatrix}\,,\quad
  {\bf d}:=\begin{pmatrix} d_1\\ \vdots\\ d_{\red{N}}\end{pmatrix}\,.
  \]
  Hence, if there exists values of these $4\red{N}$ constants
  for which the equation (\ref{E-solvability-2}) has a non-trivial
  solution, then $\Xi$ is the next vector in the generalised Jordan chain.
  
  A general theory considering all possible Jordan chains emanating from the
  condition (\ref{E-solvability-2})
  is outside the scope of this paper.  However, we will highlight special
  cases that appear in the nonlinear modulation theory.  The case
  ${\bf a}={\bf c}=0$ (a pure ${\bf J}-$chain) appears in the
  nonlinear modulation theory associated with zero characteristics
  \cite{rb16a,r17b,r18},
  and the case ${\bf b}={\bf d}=0$ is mathematically equivalent, and generates
  a pure ${\bf M}-$chain.
  \red{Here two new cases which intertwine the ${\bf J}$ and ${\bf M}$ chains,
  and are required for the nonlinear modulation theory in
  this paper, will be highlighted.}

  \red{
  \subsection{The key Jordan chain of length four}
  \label{subsec-4-2-jordanchains}}
  
Taking
\begin{equation}\label{mixed-42-a}
    {\bf a} = c^2 \be\,,\quad  
      {\bf b}={\bf c} = c \be \qand {\bf d}= \be\,,
\end{equation}
reduces the solvability condition (\ref{E-solvability-2}) to
      \begin{equation}\label{L-q-solvability}
      \left[ c^2\D_{\bw}{\bf A} + c (\D_{\bk}{\bf A} + 
  \D_{\bw}{\bf B}) + \D_{\bk}{\bf B} \right]\be = {\bf 0}\,.
      \end{equation}
      Remarkably, this is precisely the equation for characteristics.
      In this case,
      the Jordan chain associated with ${\bf L}$ can continue when $\Delta(c)=0$
      and $\be\in{\rm Ker}[{\bf E}(c)]$,
      the familiar condition (\ref{Ec-def})
      for the existence of a characteristic $c$. However, this construction
      does not imply
      that $c=c_g$, that equivalence will follow from another Jordan
      chain and it is considered in \S\ref{subsec-cg-jordanchain} below.

      In the case (\ref{mixed-42-a}) with (\ref{L-q-solvability})
      the Jordan chain intertwines
      the symplectic ${\bf J}-$chain and the symplectic
      ${\bf M}-$chain.  They can be combined to a new symplectic Jordan
      chain, based on the combined symplectic operator
      ${\bf J} + c{\bf M}$ and ultimately leads to the central
      Jordan chain that shows up in the nonlinear modulation theory.
           
      Suppose first that $c$ is arbitrary, and see that the condition
      $c=c_g$ will arise as a condition to extend the Jordan chain in
      \S\ref{subsec-cg-jordanchain} below.
      For arbitrary $c$, there is still a geometric eigenvector $\be$
      satisfying ${\bf E}(c)\be=0$.  Express it in components,
      $\be = (\zeta_1,\ldots,\zeta_{\red{N}})$, and re-number the generalized eigenvectors
      as follows,
\begin{equation}\label{v-chain-def}
\begin{array}{rcl}
  {\bf v}_1 &=& \displaystyle \sum_{j=1}^{\red{N}}\zeta_j \partial\Zh/\partial\theta_j \\[4mm]
  {\bf v}_2 &=& \displaystyle \sum_{j=1}^{\red{N}}\zeta_j (\partial\Zh/\partial k_j
  +c\,\partial\Zh/\partial \omega_j)\,,
\end{array}
\end{equation}
These two vectors satisfy
\[
  {\bf L}{\bf v}_1 =0 \qand {\bf L}{\bf v}_2 = ({\bf J}+c{\bf M}){\bf v}_1\,.
  \]
  This Jordan chain of length continues to length three if
  the following equation is solvable,
  \[
    {\bf L}{\bf v}_3 = ({\bf J}+c{\bf M}){\bf v}_2\,.
    \]
But, the existence of the ${\bf v}_3$ term is just a reformulation of the
solvability condition (\ref{E-solvability}) in terms of the new coordinates.
To see this, write out the solvability condition for ${\bf v}_3$
\[
\lth \Zh_{\theta_j},({\bf J}+c{\bf M}){\bf v}_2\rth =0\,,\quad j=1,\ldots,\red{N}\,.
\]
Using (\ref{v-chain-def}), \red{and noting that
\[
\begin{array}{rcl}
  \lth \Zh_{\theta_j},({\bf J}+c{\bf M}){\bf v}_2\rth
  &=& -\lth ({\bf J}+c{\bf M})\Zh_{\theta_j},{\bf v}_2\rth\\[2mm]
  &=& -\lth {\bf L}(\partial_{k_j}\Zh +c\partial_{\omega_j}\Zh),{\bf v}_2\rth\\[2mm]
  &=& -\lth (\partial_{k_j}\Zh +c\partial_{\omega_j}\Zh),{\bf L}{\bf v}_2\rth\\[2mm]
  &=& -\lth (\partial_{k_j}\Zh +c\partial_{\omega_j}\Zh),({\bf J}+c{\bf M}){\bf v}_1\rth\,,\quad j=1,\ldots,N\,.
\end{array}
\]
Substituting for ${\bf v}_1$ and using
the identities (\ref{Aw})-(\ref{Bk})},
generates precisely
the equation (\ref{L-q-solvability}).
Indeed it was working
backwards from (\ref{L-q-solvability}) that suggested the
definitions (\ref{v-chain-def}).  Since ${\bf L}$ is symmetric
and ${\bf J}+c{\bf M}$ is skew-symmetric every Jordan chain has even
length, assuring the existence of ${\bf v}_4$,
\[
    {\bf L}{\bf v}_4 = ({\bf J}+c{\bf M}){\bf v}_3\,.
\]
It is assumed that the this four chain terminates; that is,
the system
\begin{equation}\label{Non-solvable-chain}
  {\bf L}{\bf v}_5 = ({\bf J}+c{\bf M}){\bf v}_4\,,
  \end{equation}
  is \emph{not solvable}.
  The four chain
  \begin{equation}\label{four-chain}
    {\bf L}{\bf v}_j = ({\bf J}+c{\bf M}){\bf v}_{j-1}\,,\quad j=1,\ldots,4\,,
  \end{equation}
  with ${\bf v}_0=0$ is the Jordan chain that plays a key role in
  the nonlinear modulation theory.  The non-solvability of
  (\ref{Non-solvable-chain})
  also arises in the nonlinear modulation theory.  It ensures that
  the coefficient of dispersion, $\mathscr{K}$, is nonzero.

  In addition there are $\red{N-1}$ Jordan blocks of length two, but explicit
  expressions for these blocks are not needed in the nonlinear modulation
  theory.  It is however assumed that they are each of length exactly two.
  \vspace{.15cm}

\subsection{Another mixed Jordan chain defining $c_g$}
\label{subsec-cg-jordanchain}

There is yet another Jordan chain, associated with ${\bf L}$,
that arises in the nonlinear modulation theory and the solvability
condition for this chain defines $c_g$.  It is
a special case of the solvability condition (\ref{E-solvability-2}) obtained
by taking
\begin{equation}\label{abcd-def}
  {\bf a} = c^2\bgam + 2c\be\,, \quad
  {\bf b}  =   {\bf c}  =   c\bgam + \be\,, \quad
  {\bf d}  =  \bgam\,.
\end{equation}
Substitution of (\ref{abcd-def}) into (\ref{E-solvability-2}) and rearranging gives
  \begin{equation}\label{L-solvability-secondjordanchain}
  \left[ \D_{\bw}{\bf A}c^2 +(\D_{\bk}{\bf A}+
    \D_{\bw}{\bf B})c + \D_{\bk}{\bf B} \right]\bgam + \left[ 2c
    \D_{\bw}{\bf A}+(\D_{\bk}{\bf A}+\D_{\bw}{\bf B})\right]\be = 0\,,
  \end{equation}
  or
  \begin{equation}\label{E-chain-sec4}
    {\bf E}(c)\bgam + {\bf E}'(c)\be = 0 \,.
  \end{equation}
  This equation is solvable for a fixed value of $c$ only, and the
  solvability condition
  \[
  \langle\be,{\bf E}'(c_g)\be\rangle=0\,,
  \]
  agrees with the definition of $c_g$ in (\ref{E-chain-two}) and (\ref{cg-def}).
When $c=c_g$, the vectors
  $(\be,\bgam)$ form a Jordan chain for ${\bf E}(c_g)$ of length two.  
  
  Suppose the solvability condition (\ref{L-solvability-secondjordanchain})
  and (\ref{E-chain-sec4})
  is satisfied, then substitution back into (\ref{L.1}) gives that
\begin{equation}\label{L-xi}
{\bf L}\Xi = \sum_{i=1}^{\red{N}}\left[\zeta_i\big({\bf J}\Zh_{\omega_i}+{\bf M}\Zh_{k_i}+2c_g{\bf M}\Zh_{\omega_i}\big)+\gamma_i({\bf J}+c_g{\bf M})(\Zh_{k_i}+c_g\Zh_{\omega_i})\right]\,.
\end{equation}
It is this equation that arose in the modulation theory at fourth order
(\ref{Xi-def-eqn}) and working
backwards we see that it is a special case of (\ref{L.1}) and moreover
solvability, with the expressions (\ref{abcd-def}), is precisely the condition
for the Jordan chain (\ref{E-chain-sec4}) of ${\bf E}(c_g)$.
  
Further still, we can define another special case, which results in the criterion for the termination of this chain. This is achieved by setting
\begin{equation}\label{abcd-def_2}
  {\bf a} = c_g^2\Ups+2 c_g\bgam + \be\,, \quad
  {\bf b}  =   {\bf c}  = c_g\Ups+\bgam \,, \quad
  {\bf d} = \Ups\,.
\end{equation}
Utilising this in (\ref{E-solvability-2}) and simplifying results in the system
\[
    {\bf E}(c_g)\Ups+{\bf E}'(c_g)\bgam + \fr{\bf E}''(c_g)\be = 0\,,
\]
which is precisely
 (\ref{Upsilon-eqn}). The assumption made here is that this chain is of length two, 
 and so the right hand side of (\ref{Upsilon-eqn}) does not lie in the range of ${\bf L}$. 
 Thus, by appealing to solvability one recovers the condition that
 \[
\mu:=\langle\be,{\bf E}'(c_g)\bgam + \fr{\bf E}''(c_g)\be\rangle =\fr \langle\be,{\bf E}''(c_g)\be\rangle - \langle\bgam,{\bf E}(c_g)\bgam \rangle\neq 0\,,
 \]
 and therefore completing the connection between $\mu$ and the termination of this
 mixed Jordan chain. Within the modulation theory, this corresponds to the system
 \[
 {\bf L}\Gamma =  \sum_{i=1}^{\red{N}}\left[\zeta_i{\bf M}\Zh_{\omega_i}+\gamma_i\big({\bf J}\Zh_{\omega_i}+{\bf M}\Zh_{k_i}+2c_g{\bf M}\Zh_{\omega_i}\big)+\Upsilon_i({\bf J}+c_g{\bf M})(\Zh_{k_i}+c_g\Zh_{\omega_i})\right]\,,
 \]
being unsolvable for $\Gamma$, and what ultimately leads to the coefficient of the time derivative term in the emergent two-way Boussinesq equation.
 
We have only scratched the surface of the possible solvability conditions
and attendant Jordan chains associated with (\ref{L.1}).  However, we 
have all the Jordan chains needed for the nonlinear modulation theory.

\section{Properties of the two-way Boussinesq equation}
\setcounter{equation}{0}
\label{sec-boussinesq-eqn}

Once the modulation equation
is derived in a specific context, analysis of the solutions of
the two-way Boussinesq equation (\ref{U-2way-Boussinesq})
gives information about the nature of solutions in the
nonlinear problem near coalescence.

The two-way Boussinesq equation is valid at $c=c_g$. At least one
parameter needs to be varied to obtain the coalescence.  That
parameter can be a one-parameter path through the
four-dimensional frequency-wavenumber
$(\bw,\bk)$ space, or it could be a perturbation of the frame speed
$c=c_g+\mathcal{O}(\eps^2)$.  Unfolding the singularity
generates a term of the form $\nu U_{XX}$ in (\ref{U-2way-Boussinesq}),
regardless of the precise perturbation
path (this can be shown by perturbing
the linearised generic Whitham equations).  Therefore the full
modulation equation in the neighbourhood of the coalescence is
\begin{equation}\label{U-2way-Boussinesq-nu}
  \mu U_{TT}  + \nu U_{XX} + \kappa\, (UU_X)_X + \mathscr{K} U_{XXXX}=0\,\red{,}
\end{equation}
where $\nu$ is an order one constant.

\begin{figure}[htbp]
  \centering
\includegraphics[width=11.0cm]{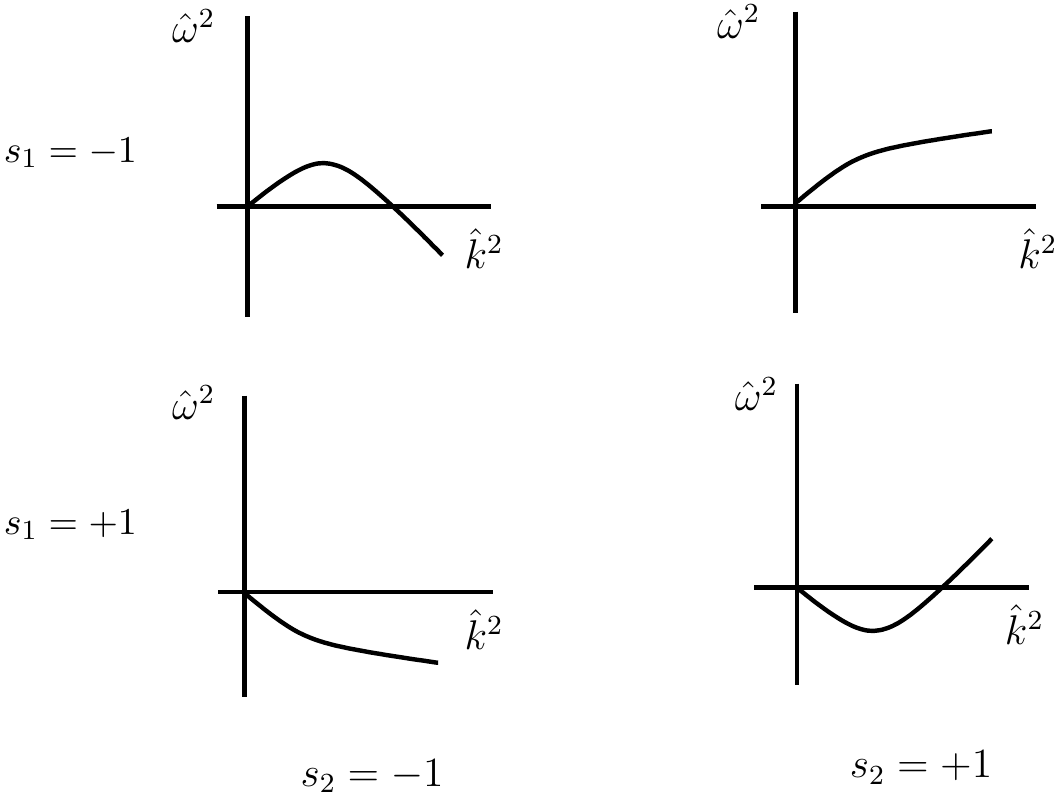}
\caption{The four qualitative cases of the dispersion relation
  (\ref{twoway-B-disp-rel}) determined by the signs $s_1$
  and $s_2$ in the two-way Boussinesq equation (\ref{two-way-boussinesq-app}).}
\label{fig-LPfourcases}
\end{figure}

When the coefficients are non-zero, the Boussinesq equation can be put into
standard form.  Scale the independent and dependent variables:
$\tau=aT$, $\xi=b X$, and $U = \rho u$; then values of $a,b,\rho$ can be
chosen so that the two-way Boussinesq equation becomes
\begin{equation}\label{two-way-boussinesq-app}
  u_{\tau\tau} +  s_1 u_{\xi\xi} + \big(\fr u^2)_{\xi\xi} + s_2 u_{\xi\xi\xi\xi} = 0 \,,
  \quad s_1,s_2=\pm 1\,,
\end{equation}
with
\[
s_1 = {\rm sign}(\mu\nu) \qand s_2 = {\rm sign}\left(\mu\mathscr{K}\right)\,.
\]
The sign $s_1$ determines whether the unfolding is into the
elliptic region ($\red{s_1}=+1$) or into the
hyperbolic region ($\red{s_1}=-1$, in which case all characteristics are hyperbolic).
The sign $s_2$ indicates whether the resulting
two-way Boussinesq equation is good ($s_2=+1$) or bad ($s_2=-1$).
In the latter case, the initial value problem for the linearized
system about $u=0$ is \red{ill-posed}, and small initial data
 with zero mean is therefore expected to saturate to form
 nonlinear structures.  The ill-posedness in the case $s_2=-1$
 can be seen by
considering the linearization of (\ref{two-way-boussinesq-app}) about
the trivial solution and introducing a normal mode solution of the form
$\re^{\ri(\hat k\xi+\hat\omega\tau)}$.  The dispersion relation
associated with the normal mode is then
\begin{equation}\label{twoway-B-disp-rel}
\hat\omega^2 = -s_1 \hat k^2 + s_2 \hat k^4\,.
\end{equation}
There are four cases depending on the signs $s_1$ and $s_2$, and they
are shown in Figure \ref{fig-LPfourcases}.  The figure plots
$\hat\omega^2$ against $\hat k^2$ and so $\hat\omega^2<0$ indicates
linear instability of the trivial solution which in turn reflects
linear instability of the basic travelling wave.

When $s_1<0$ (the upper two cases in Figure \ref{fig-LPfourcases})
then either an unstable band
emerges at finite $\hat k$ when $s_2=-1$ or the Boussinesq
equation is 
hyperbolic for all wavenumbers ($s_2=+1$).
When $s_1>0$ (lower two cases in Figure
\ref{fig-LPfourcases}) then
either a cutoff wave number emerges with re-stabilization at finite $\hat k$
(as in the lower right diagram with $s_2=+1$), or instability is further
enhanced for all wavenumbers ($s_1=+1$ and $s_2=-1$).

The simplest class of nonlinear solutions of (\ref{two-way-boussinesq-app}) are
travelling solitary wave solutions, for example,
\[
u(\xi,\tau) = \widehat u(\xi+\gamma\tau)\,,
\]
which satisfies the ODE
\[
\big( \gamma^2 \widehat u + s_1 \widehat u + \fr \widehat u^2
+ s_2 \widehat u''\big)'' = 0\,.
\]
Integrating and taking the function of integration to be constant
\[
s_2\widehat u'' + (s_1+\gamma^2)\widehat u + \fr\widehat u^2 = h\,.
\]
The constant of integration $h$ is fixed by initial data or
the value of $\widehat u$ at infinity.
For appropriate parameter values,
this planar ODE has a family of periodic solutions
and a homoclinic orbit which represent periodic travelling
waves and a solitary travelling wave solution of
(\ref{two-way-boussinesq-app}).
The implication of these solutions is that the
transition from elliptic to hyperbolic of a periodic travelling
wave of the original system generates a coherent structure in the
transition, which is represented by the above solitary wave.  However,
there is much more complexity generated at the transition.
\textsc{Hirota}~\cite{hirota} shows that there is a large family of
$M-$soliton solutions to (\ref{two-way-boussinesq-app}) as well,
where $M$ can be any natural number.
Further details especially in the case $M=2$ are given in \cite{hirota}.
Blow-up can occur in the two-way Boussinesq equation even in the
case of the good Boussinesq equation \cite{turitsyn}.
It is also generated by a Lagrangian,
and has both a Hamiltonian and multisymplectic structure
(e.g.\ \S10 of \cite{bd01} and \cite{chen}).

\section{CNLS wavetrains with coalescing characteristics}
    \setcounter{equation}{0}
    \label{sec-coupled-nls}

    To illustrate the \red{nonlinear} theory it is applied to the modulation of
    two-phase wavetrains of a coupled nonlinear Schr\"odinger (CNLS) equation.
    This example serves two purposes: firstly,
    it shows that the coalescence of characteristics
    is quite common and appears even in the simplest of examples,
    and secondly, it shows that computing the coefficients in the
    emergent two-way Boussinesq equation is elementary once the
    properties of the basic state are found.
    
    The CNLS equation is a
    canonical example of a PDE generated by a Lagrangian with a
    toral symmetry, $\T^2=S^1\times S^1$.  Indeed any finite number of
    NLS equations can be coupled together to generate
    a toral symmetry $\T^N$ for any natural
    number $N$, and they will have explicit $N-$phase wavetrains
    which are also relative equilibria.
    Here attention is restricted to two coupled NLS equations in the form
  \begin{equation}\label{CNLS}
  \begin{array}{rcl}
\displaystyle    \ri\frac{\partial \Psi_1}{\partial t} +
    \alpha_1 \frac{\partial^2 \Psi_1}{\partial x^2} +
    (\beta_{11}|\Psi_1|^2 + \beta_{12}|\Psi_2|^2)\Psi_1 &=& 0 \\[4mm]
    \displaystyle    \ri\frac{\partial \Psi_2}{\partial t} +
    \alpha_2 \frac{\partial^2 \Psi_2}{\partial x^2} +
    (\beta_{21}|\Psi_1|^2 + \beta_{22}|\Psi_2|^2)\Psi_2 &=& 0\,,
  \end{array}
  \end{equation}
  where the coefficients $\alpha_j,\beta_{ij}$, $i,j=1,2$, are given
  real constants, with $\beta_{21}=\beta_{12}$.  The
  functions $\Psi_j(x,t)$ are complex-valued and $\ri^2=-1$.

  Coupled NLS equations appear in a wide range of applications.  Two
  applications that motivated this work are the coupled NLS equations
  that appear in the theory of water waves
  (e.g.\ \textsc{Roskes}~\cite{roskes},
  \textsc{Ablowitz \& Horikis}~\cite{ah15}, \textsc{Degasperis et al.}~\cite{dls19}), and in models for Bose-Einstein condensates (e.g.\ \textsc{Salman \& Berloff}~\cite{sb09}, \textsc{Kevrekidis \& Frantzeskakis}~\cite{kf16}).
    The PDE (\ref{CNLS}) is the Euler-Lagrange equation for
  \[
  \mathcal{L}(\Psi) = \int_{t_1}^{t_2}\int_{x_1}^{x_2} L(\Psi_t,
  \Psi_x,\Psi)\,\rd x \rd t\,,
  \]
  with $\Psi:=(\Psi_1,\Psi_2)$ and
  \[
  \begin{array}{rcl}
    L &=& \displaystyle\frac{\ri}{2}\left(\overline\Psi_1(\Psi_1)_t -
    \Psi_1(\overline\Psi_1)_t \right) +
    \frac{\ri}{2}\left(\overline\Psi_2(\Psi_2)_t -
    \Psi_2(\overline\Psi_2)_t \right) \\[4mm]
    &&\quad -\alpha_1\big|(\Psi_1)_x\big|^2
    -\alpha_2\big|(\Psi_2)_x\big|^2 +\fr\beta_{11}|\Psi_1|^4 +
    \beta_{12}|\Psi_1|^2|\Psi_2|^2 + \fr\beta_{22}|\Psi_2|^4\,,
  \end{array}
  \]
  with the overline denoting complex conjugate.

  The toral symmetry follows from the fact that
  $(\re^{\ri\theta_1}\Psi_1,\re^{\ri\theta_2}\Psi_2)$ is a solution of
  (\ref{CNLS}), for any $(\theta_1,\theta_2)\in S^1\times S^1$,
  when $(\Psi_1,\Psi_2)$ is a solution.  The complex coordinates can
  be converted to real coordinates, generating a standard action of
  $\T^2$ but will not be needed as the main calculations can be
  done in the complex setting.
  
  Noether's theorem gives the conservation laws
  \begin{equation}\label{AB-nls-claw}
  (A_j)_t + (B_j)_x =0\,,\quad j=1,2\,,
  \end{equation}
  with
  \begin{equation}\label{AB-nls-components}
  A_j = \fr|\Psi_j|^2 \qand B_j = \alpha_1{\rm Im}(\overline\Psi_1(\Psi_1)_x)\,,
  \quad j=1,2\,.
  \end{equation}
  The basic state is just the usual family of plane waves, but
  interpreted here as a family of relative equilibria associated with
  the $\T^2$ symmetry; it has the form,
  \begin{equation}\label{NLS-wavetrains}
  \Psi_j(x,t) = \Psi_j^0(\bw,\bk)\re^{\ri\theta_j(x,t)}\,,\quad
  \theta_j(x,t) = k_jx+\omega_jt+\theta_j^0\,,\quad j=1,2\,.
  \end{equation}
  Substitution into the governing equations (\ref{CNLS}) generates
  the required relationship between the amplitudes, frequencies and
  wavenumbers,
\begin{equation}\label{basic-state}
\begin{array}{rcl}
|\Psi_1^0|^2 &=&\displaystyle 
\frac{1}{\beta}\bigg(\beta _{22}(\omega_1+\alpha _1 k_1^2)-\beta _{12}(\omega_2+\alpha _2 k_2^2)\bigg)\\[4mm]
|\Psi_2^0|^2 &=&\displaystyle \frac{1}{\beta}\bigg(\beta _{11}(\omega_2+\alpha _2 k_2^2)-\beta _{21}(\omega_1+\alpha_1 k_1^2)\bigg)\,,
\end{array}
\end{equation}
with $\beta_{21}:=\beta_{12}$ and $\beta=
\beta _{11}\beta _{22}-\beta _{12}\beta_{21}\neq 0$.

The key wave action vectors
${\bf A}(\bw,\bk)$ and ${\bf B}(\bw,\bk)$,
needed for analysis of the linearization, are obtained by substituting (\ref{basic-state})
into the components of the conservation law (\ref{AB-nls-components}),
\begin{equation}\label{A-script}
 {\bf A}(\bw,\bk) := \begin{pmatrix}\mathscr{A}_1(\bw,\bk)\\
   \mathscr{A}_2(\bw,\bk)\end{pmatrix} =
 \frac{1}{2\beta}\begin{pmatrix}
   \beta_{22}(\omega_1+\alpha_1 k_1^2)-\beta_{12}(\omega_2+\alpha _2 k_2^2)\\
   \beta_{11}(\omega_2+\alpha_2 k_2^2)-\beta_{21}(\omega_1+\alpha_1 k_1^2) \end{pmatrix}
\end{equation}
and
\begin{equation}\label{B-script}
 {\bf B}(\bw,\bk) := \begin{pmatrix}\mathscr{B}_1(\bw,\bk)\\
   \mathscr{B}_2(\bw,\bk)\end{pmatrix} =
 \frac{\alpha_1 k_1}{\beta}\begin{pmatrix}
   \beta_{22}(\omega_1+\alpha_1 k_1^2)-\beta_{12}(\omega_2+\alpha_2 k_2^2)\\
   \beta_{11}(\omega_2+\alpha_2 k_2^2)-\beta_{21}(\omega_1+\alpha_1 k_1^2)
 \end{pmatrix}\,.
\end{equation}
The linear operator ${\bf E}(c)$ defined in (\ref{wme-3-quadratic})
is 
\[
{\bf E}(c) := \D_{\bw}{\bf A}c^2 + (\D_{\bk}{\bf A}+\D_{\bw}{\bf B})c +
  \D_{\bk}{\bf B}\,,
\]
  with
  \begin{equation}\label{Dw-def}
\D_{\bw}{\bf A} = \frac{1}{2\beta}\begin{pmatrix} \beta_{22} & -\beta_{12}\\
  -\beta_{12} & \beta_{11} \end{pmatrix} \,,
  \end{equation}
and
  \begin{equation}\label{Dk-def}
\D_{\bk}{\bf A} = \frac{1}{\beta}\begin{pmatrix}
  \alpha_1\beta_{22}k_1 & -\alpha_2\beta_{12}k_2 \\
  -\alpha_1\beta_{12}k_1 & \alpha_2\beta_{11}k_2 \end{pmatrix} = \D_{\bw}{\bf B}^T\,,
  \end{equation}
and
\begin{equation}\label{D-k-B}
{\rm D}_{\bf k}{\bf B} = \frac{1}{\beta}
\begin{pmatrix}\alpha_1\beta|\Psi_1^0|^2+2\beta_{22}\alpha _1^2 k_1^2 &
 -2\beta_{12}\alpha _1 \alpha _2 k_1 k_2 \\[2mm]
-2\beta _{12}\alpha _1 \alpha _2 k_1 k_2 &
\alpha_2\beta|\Psi_2^0|^2 + 2\alpha _2^2 \beta _{11}k_2^2
\end{pmatrix}\,.
\end{equation}
The characteristic polynomial is
\begin{equation}\label{nls-char-poly}
  \Delta(c) := {\rm det}[{\bf E}(c)] = a_0c^4 + a_1c^3 + a_2 c^2 + a_3 c + a_4\,,
\end{equation}
with
\begin{equation}\label{cnls-quartic}
\begin{array}{rcl}
a_0 &=& \frac{1}{4}\beta^{-1}\,,\\[3mm]
a_1 &=& \beta^{-1}(\alpha_1k_1+\alpha_2k_2)\,,\\[3mm]
a_2&=& \frac{1}{2}\beta^{-1}\big[\alpha_1(\beta_{11}|\Psi_1^0|^2+2\alpha_1k_1^2)+\alpha_2(\beta_{22}|\Psi_2^0|^2+2\alpha_2k_2^2)+8\alpha_1\alpha_2k_1k_2\big]\,,\\[3mm]
a_3 &=& 2\alpha_1\alpha_2 \beta^{-1}\big(k_1(\beta_{22}|\Psi_2^0|^2+2\alpha_2k_2^2)+k_2(\beta_{11}|\Psi_1^0|^2+2\alpha_1k_1^2)\big)\\[3mm]
a_4 &=& \alpha_1\alpha_2\beta^{-1}\big((\beta_{11}|\Psi_1^0|^2+2\alpha_1k_1^2)(\beta_{22}|\Psi_2^0|^2+2\alpha_2k_2^2)-|\Psi_1^0|^2|\Psi_2^0|^2\beta_{12}^2\big)\,.
\end{array}
\end{equation}
Coalescing characteristics are obtained by solving $\Delta(c)=\Delta'(c)=0$
for $c$.  This problem is solved numerically in \cite{br19} by using \red{the}
\emph{graphical sign characteristic}.  The function $\Delta(c)$ is
plotted versus $c$ as parameters vary.  That way roots and points
where $\Delta'(c)=0$ can be read off the graph.  It is inspired by the
graphical Krein signature introduced by \textsc{Kollar \& Miller}~\cite{km14}.
Results in \cite{br19} show that coalescing characteristics are plentiful
in the Whitham modulation theory for CNLS.  

According to the theory in this paper at coalescing characteristics
the following nonlinear modulation equation is generated
\begin{equation}\label{two-way-B-CNLS}
\mu U_{TT} + \kappa(UU_X)_X + \mathscr{K}U_{XXXX}=0\,.
\end{equation}
In principle the quartic $\Delta(c)=0$ can be solved in closed form,
but in practice this is lengthy and not illuminating, and numerical
methods are more effective.  For simplicity here, the
case of standing waves (where $\bk=0$) are considered,
which restricts the parameter space
significantly, and so calculations can be done explicitly.
The strategy for calculating  $\mu$ and
$\kappa$ is to construct the averaged Lagrangian and use the formulas
(\ref{kappa-def}) and (\ref{alpha-def}).

\subsection{Calculations for standing waves}
\label{subsec-standing-waves}

Standing waves are defined as basic states of the form (\ref{NLS-wavetrains})
but with $\bk=0$.  With this restriction the coefficients $a_1$ and $a_3$
are identically zero reducing the coefficients in the
polynomial in (\ref{cnls-quartic}) to
\[
\begin{array}{rcl}
a_0 &=& \frac{1}{4}\beta^{-1}\,,\\[3mm]
a_2&=& \frac{1}{2}\beta^{-1}\big[\alpha_1\beta_{11}|\Psi_1^0|^2+\alpha_2\beta_{22}|\Psi_2^0|^2\big]\,, \\[3mm]
a_4 &=&\alpha_1\alpha_2|\Psi_1^0|^2|\Psi_2^0|^2\,.
\end{array}
\]
There are four characteristics and they satisfy the biquadratic equation
\[
a_0 c^4 + a_2 c^2 + a_4 = 0\,,
\]
giving
\begin{equation}\label{cg-squared}
c^2 =  -\alpha_1\beta_{11}|\Psi_1^0|^2-\alpha_2\beta_{22}|\Psi_2^0|^2 \pm \sqrt{(\alpha_1\beta_{11}|\Psi_1^0|^2-\alpha_2\beta_{22}|\Psi_2^0|^2)^2+4\alpha_1\alpha_2\beta_{12}^2|\Psi_1^0|^2|\Psi_2^0|^2} \,.
\end{equation}
Coalescing characteristics occur precisely when
the discriminant vanishes
\[
(\alpha_1\beta_{11}|\Psi_1^0|^2-\alpha_2\beta_{22}|\Psi_2^0|^2)^2+4\alpha_1\alpha_2\beta_{12}^2|\Psi_1^0|^2|\Psi_2^0|^2= 0\,.
\]
One way to interpret this equation is as
a line in the positive quadrant of
$\big( |\Psi_1^0|^2 ,|\Psi_2^0|^2\big)$ space defined by
\begin{equation}\label{cc-line}
\alpha_2\beta_{22}^2|\Psi_2^0|^2 = \alpha_1\big(\beta_{11}\beta_{22}-2 \beta_{12}^2 \pm 2 \red{|\beta_{12}|}\sqrt{-\beta\,}\,\big) |\Psi_1^0|^2\,,
\end{equation}
which includes the conditions $\beta<0$ and $\alpha_1\alpha_2<0$ for
reality.  At coalescence it follows from (\ref{cg-squared}) that
\[
c_g^2 = -\alpha_1\beta_{11}|\Psi_1^0|^2 -\alpha_2\beta_{22}|\Psi_2^0|^2\,,
\]
which carries with it the requirement that
$\alpha_1\beta_{11}|\Psi_1^0|^2 +\alpha_2\beta_{22}|\Psi_2^0|^2<0$,
a condition that is effectively a generalisation of the defocussing classification for the one-component NLS.

Now suppose parameters are such that (\ref{cc-line}) is satisfied,
and proceed to compute the required coefficients in (\ref{two-way-B-CNLS}).
The eigenvector and generalised eigenvector of ${\bf E}(c_g)$ are,
\begin{equation}
\begin{split}
\be &=
\begin{pmatrix}
c_g^2 \beta_{12}\\
\beta_{22}c_g^2+2 \alpha_1\beta |\Psi_1^0|^2 
\end{pmatrix}\,, \\[3mm] 
\bgam &= -\frac{8c_g\alpha_1 \beta_{12}|\Psi_1^0|^2\beta^2}{\beta_{22}c_g^2+2\alpha_1\beta|\Psi_1^0|^2}
\begin{pmatrix}
1\\
0
\end{pmatrix}\,.
\end{split}
\end{equation}
Now use these eigenvectors and the Jacobians (\ref{Dw-def}), (\ref{Dk-def}) and
(\ref{D-k-B}) to generate the coefficients of the emergent Boussinesq equation.
The first computed is the coefficient of the time derivative term,
\[
\be^T\D_\bw{\bf A}\be+\be^T{\bf E}'(c_g)\bgam = 4 c_g^2 \kappa_0\,, \quad {\rm with} \quad  \kappa_0 = 4 \beta (\beta_{22}c_g^2+2 \alpha_1\beta |\Psi_1^0|^2 )\,.
\]
Next, one may use the variation of the Lagrangian to show that the coefficient of the nonlinear term is
\[
\kappa =-\frac{3c_g^2\kappa_0}{2 |\Psi_2^0|^2}(\alpha_1\beta_{11}|\Psi_1^0|^2-\alpha_2 \beta_{22}|\Psi_2^0|^2)(\alpha_1\beta_{11}|\Psi_1^0|^2-\alpha_2 \beta_{22}|\Psi_2^0|^2+2 \alpha_1 \beta_{12}|\Psi_2^0|^2)\,,
\]
The coefficient of dispersion requires a Jordan chain analysis. This
would require multisymplectification of CNLS and construction of the
linear operator ${\bf L}$.  However, this CNLS has been multisymplectified
in \textsc{Ratliff}~\cite{r17b}, where reduction to
KdV and 2-parameter Boussinesq were studied, and the Jordan chain
theory is close to this case.
With minor modification of that analysis,
the desired dispersive coefficient is found to be
\[
\be^T{\bf T}= \frac{\kappa_0(\alpha_2 |\Psi_1^0|^2 -\alpha_1|\Psi_2^0|^2) }{|\Psi_1^0|^2 |\Psi_2^0|^2 (\alpha_1\beta_{11}|\Psi_1^0|^2 - \alpha_2\beta_{22}|\Psi_2^0|^2)}\,.
\]
Each of these coefficients has a common factor $\kappa_0$, and so the
two-way Boussinesq that emerges at the coalescence of characteristics
simplifies to
\begin{equation}\label{2wayB-cnls}
c_g^2U_{TT}+\bigg(\fr\tilde\kappa U^2+\tilde{\mathscr{K}} U_{XX}\bigg)_{XX} = 0\,,
\end{equation}
with
\[
\begin{split}
\tilde\kappa &= -\frac{3c_g^2}{8|\Psi_2|^2}(\alpha_1\beta_{11}|\Psi_1|^2-\alpha_2 \beta_{22}|\Psi_2|^2)(\alpha_1\beta_{11}|\Psi_1|^2-\alpha_2 \beta_{22}|\Psi_2|^2+2 \alpha_1 \beta_{12}|\Psi_2|^2)\,,\\[2mm]
\tilde{\mathscr{K}}&=\frac{\alpha_2 |\Psi_1^0|^2 -\alpha_1|\Psi_2^0|^2}{4|\Psi_1^0|^2 |\Psi_2^0|^2 (\alpha_1\beta_{11}|\Psi_1^0|^2 - \alpha_2\beta_{22}|\Psi_2^0|^2)}\,.
\end{split}
\]
With $\tilde\kappa$ and $\tilde{\mathscr{K}}$ nonzero,
one can proceed to analyze the solutions of this equation using
results in the literature (e.g.\ \cite{hirota,turitsyn}).
A detailed analysis of (\ref{2wayB-cnls}) and its implications for
coupled NLS is outside the scope of this paper, but the diversity
of complexity due to coalescing characteristics is
clear; for example, evaluation of $\tilde{\mathscr{K}}$ along the
lines (\ref{cc-line}) shows that (\ref{2wayB-cnls}) 
can be both positive (good Boussinesq) and negative (bad Boussinesq).

\section{Concluding remarks}
\setcounter{equation}{0}
\label{sec-cr}

This paper gives a complete weakly nonlinear theory for multiphase
WMT when a pair of characteristics coalesce and transition from
hyperbolic to elliptic.  This transition, in the nonlinear problem
creates  nonlinear dispersive dynamics, and it transpires that the resulting normal form is the two-way Boussinesq equation.  \red{There are 
potential generalisations and new directions emerging from this theory.}

\subsection{Generalisation to $2+1$}
\label{cr-twoplusone}

Although we have confined the discussion to $1+1$ dimensions,
 there is a natural generalisation to $2+1$.  A good starting point
 is the $2+1$ theory for the nonlinear modulation of
 single-phase wave trains near coalescing characteristics
 \cite{br18}.  However, the Jordan chain theory
 in (\ref{L.1}) will literally take on
 a new dimension, bringing in the intertwinement of three symplectic
 Jordan chains.  On the other hand, key features like the frame speed,
 scaling, sign characteristic,
 and reduction should carry over with appropriate modification.

  \subsection{Examples}
 \label{cr-examples}
 
 The results in \red{this} paper are universal, and are operational whenever 
  a Lagrangian system has a suitable characteristic collision, which can be identified
  via the sign characteristic diagnostic used in \cite{br19} and
 in this paper.
There are examples in the literature where multiphase Whitham modulation
theory has been applied and coalescing characteristics observed, and
so the application of the theory in this paper is relevant.  Two
examples are Stokes travelling waves coupled to meanflow
(\textsc{Whitham}~\cite{whitham67}, \textsc{Willebrand}~\cite{willebrand})
and modulation of viscous
conduit periodic waves (\textsc{Maiden \& Hoefer}~\cite{mh16}).
Both of these examples have special features which require
additional methodology.  In the case of viscous conduit waves \cite{mh16}
the equations are not generated by a Lagrangian so
the theory would have to be built on averaging of conservation laws.
However at coalescing characteristics
one expects a two-way Boussinesq equation to be generated or analogous equation
with additional non-conservative terms.
The case of modulation of Stokes waves in
shallow water \cite{whitham67} involves the full water wave problem and so the
class of PDEs (\ref{MJS}) has to be modified to account for the
vertical variation of water wave fields.  However the full water wave problem
has a multisymplectic structure (e.g.\ Chapter 14 of \cite{tjb-book}) and
so the theory should go through as in this paper, with appropriate
modification.  

\subsection{Larger kernel of ${\bf E}(c_g)$}
\label{cr-Ekernel}

\red{In this paper the basic state has $N-$phases but the dimension of the kernel of
the $N\times N$ matrix ${\bf E}(c_g)$ is one.}
A different problem arises when the kernel of
${\bf E}(c_g)$ has dimension greater
than one. In this case, the secondary reduction to ${\rm span}\{\be\}$ would
be modified to ${\rm span}\{\be_1,\ldots,\be_k\}$ where $k\leq N$ is the dimension of
the kernel of ${\bf E}(c_g)$.  Then $k-$additional coupled modulation equations
are generated (one linked to each kernel direction).

\subsection{Moving frames}
 \label{cr-movingframes}
 
 Whitham theory can also be formulated relative to any moving frame,
 and some frame speeds are more interesting than others.
 \textsc{Ratliff}~\cite{r19} shows that
 even generic Whitham theory, in the hyperbolic
 case, re-modulated relative to the appropriate characteristic frame,
 generates dispersion, on a longer time scale.

\subsection{Higher order singularities}
\label{cr-higherordersingularities}

Even in the case of two phases the parameter space is at least four dimensional,
involving $\omega_1,\omega_2,k_1,k_2$, with further degrees of freedom emerging 
when system parameters are present.  Hence higher order singularities
are to be expected, e.g. more than two characteristics coalescing, or the
coefficients $\mu$, $\kappa$, and $\mathscr{K}$ passing through zero.
A potential rescaling and re-modulation could then be implemented
leading to (as yet unknown) modulation equations replacing the two-way
Boussinesq equation.

\subsection{Hyperbolicity of multiphase modulation}
\label{cr-hyperbolicity}

For quadratic Hermitian matrix pencils a general condition for hyperbolicity can
be given.  Hyperbolicity meaning all real characteristics.  Consider
the $N-$phase case for
\begin{equation}\label{E-N-def}
{\bf E}(c){\bf u}=0\quad\mbox{with}\quad
{\bf E}(c) =\D_{\bw}{\bf A} c^2 + c \left(\D_{\bw}{\bf B}+\D_{\bk}{\bf A}\right)
  c + \D_{\bk}{\bf B}\,.
\end{equation}
 Let ${\bf u}\in\C^N$ be arbitrary and
  define
  \[
  \begin{array}{rcl}
    \alpha &=& \langle{\bf u},\D_{\bw}{\bf A}{\bf u}\rangle\\[2mm]
    \beta &=& \fr\langle{\bf u},\left( \D_{\bw}{\bf B}+\D_{\bk}{\bf A}\right)
    {\bf u}\rangle\\[2mm]
    \gamma &=& \langle{\bf u},\D_{\bk}{\bf B}{\bf u}\rangle\,,
  \end{array}
  \]
  with $\langle\cdot,\cdot\rangle$ an inner product on $\C^N$.
  \textsc{Guo \& Lancaster}~\cite{gl05} study quadratic eigenvalue
  problems in general and applying their definition to (\ref{E-N-def})
  gives the following.
  \vspace{.15cm}

  \noindent{\bf Definition}. {\it The quadratic Hermitian matrix pencil
    (\ref{E-N-def}) is hyperbolic if
  $\beta^2 > \alpha\gamma$ for all nonzero ${\bf u}\in\C^n$.}
  \vspace{.15cm}
  
  \noindent If this condition is satisfied then all the characteristics
  are real, and no coalescence can occur.  It is expected
  that the absence of coalescence would be rare.  The CNLS example shows
  coalescence to be quite common, already with $N=2$.  For arbitrary $N$,
  the parameter space $(\bw,\bk)$ has dimension $2N$ and so there is a high
  probability of coalescence.
  On the other hand, the above definition is a useful
  starting point in the analysis of multiphase WMEs.
  In the paper \cite{gl05} they go on to give a number of sufficient
  conditions, and an
  algorithm for testing hyperbolicity and computing all the eigenvalues.
These algorithms may be helpful in the study of the characteristics
of multiphase Whitham theory.

  \subsection{Number of phases tending to infinity}
  \label{cr-infinitephases}

  There is a known case where multiphase Whitham equations are hyperbolic. 
  The paper of \textsc{Willebrand}~\cite{willebrand} derives the
  multiphase WMEs and takes the limit
  $N\to\infty$ and argues that they are hyperbolic
  in this limit.  
  The argument proceeds by formally constructing explicit expressions
  for the leading order nonlinear corrections.  Small divisors
  and divergence are expected, but only the leading order terms are
  studied.  When $N$ is small, ``splitting of group velocity''
  is noted in the weakly nonlinear case, which is equivalent to
  what is called ``coalescing characteristics'' in this paper. The unfolding
  of this split group velocity may lead to instability. But
  \textsc{Willebrand} argues that the splitting disappears as $N\to\infty$.
     In the context of this paper the limit
  $N\to\infty$ would just replace the matrix pencil ${\bf E}(c)$ by an
  Hermitian operator pencil and so Willebrand's claim would be that
  ${\bf E}(c)$ in the case $N\to\infty$ is hyperbolic.
  It is important to keep in mind that this argument is
  for multiphase modulation of weakly nonlinear Stokes waves only, but
  is an intriguing example nevertheless.

\bibliographystyle{amsplain}

\end{document}